\documentclass[preprint,3p,12pt]{elsarticle}

\usepackage{stmaryrd}
\usepackage{bbding}
\usepackage{dcolumn}
\usepackage{graphicx}
\usepackage{amsmath}
\usepackage{amsfonts}
\usepackage{amssymb}
\usepackage{psfrag}
\usepackage{wrapfig}
\usepackage{subfigure}
\usepackage{makeidx}
\usepackage{bm}
\usepackage{epsf}
\usepackage{epsfig}

\begin{document}

\title{A Compact Third-order Gas-kinetic Scheme for Compressible Euler and Navier-Stokes Equations}

\author{Liang Pan}
\ead{panliangjlu@sina.com}
\author{Kun Xu\corref{cor1}}
\ead{makxu@ust.hk} \cortext[cor1]{Corresponding author}

\address{Department of mathematics, Hong Kong University of Science and Technology, Clear Water Bay, Kowloon, Hong Kong}

\begin{abstract}
In this paper, a compact third-order  gas-kinetic scheme is proposed
for the compressible Euler and Navier-Stokes equations.
The main reason for the feasibility to develop such a high-order scheme with compact stencil, which involves only neighboring cells,
is due to the use of a high-order gas evolution model. Besides the evaluation of the time-dependent flux function across a cell interface, the high-order
gas evolution model also provides an accurate time-dependent solution of the flow variables at a cell interface.
Therefore, the current scheme not only updates the cell averaged conservative flow variables
inside each control volume, but also tracks the flow variables at the cell interface at the next time level.
As a result, with both cell averaged and cell interface values the high-order reconstruction in the current scheme can be done compactly.
Different from using a weak formulation for high-order accuracy in the Discontinuous Galerkin (DG) method,
the current scheme is based on the strong solution, where the flow evolution
starting from a piecewise discontinuous high-order initial data is precisely followed.
The cell interface time-dependent flow variables can be used
for the initial data reconstruction at the beginning of next time step.
Even with compact stencil, the current scheme has third-order accuracy in the smooth flow regions, and has favorable shock capturing property in the
discontinuous regions.
We believe that the current scheme is one of the most robust and accurate third-order compact schemes for both smooth and discontinuous
viscous and heat conducting flow simulations. It can be faithfully used from the incompressible limit to the hypersonic flow computations.
Many test cases are used to validate the current scheme.
In comparison with many other high-order schemes, the current method avoids the use of Gaussian points for the flux evaluation along the cell interface and the multi-stage
Runge-Kutta time stepping technique. Even with the increasing of computational cost in the evaluation of a multidimensional
time-dependent gas distribution function at a cell interface,
the current scheme is still efficient. Also, due to its multidimensional property of including both derivatives of flow variables in the normal and
tangential directions of a cell interface, the viscous flow solution, especially those with vortex structure, can be accurately captured.
With the same stencil of a second order scheme, numerical tests clearly demonstrate that
the current compact third-order scheme is as robust as well-developed second-order shock capturing schemes, but provides much
more accurate numerical solutions than the second order counterparts.

\end{abstract}


\maketitle

\section{Introduction}
Most computational fluid dynamics methods used in
practical applications are second-order methods. They are generally
robust and reliable. For the same computing cost, the high-order
methods (order $\geq 3$) can provide more accurate solutions, but they are less
robust and more complicated. There has been a surge of research
activities on the development of high-order methods for solving the Euler
and Navier-Stokes equations. At the current stage, many high-order
numerical methods have been developed, including discontinuous
Galerkin (DG), spectral volume (SV), spectral difference (SD),
correction procedure using reconstruction (CPR), essential
non-oscillatory (ENO), and weighted essential non-oscillatory (WENO),
etc. The DG scheme was first proposed in \cite{DG1} to solve the
neutron transport equation. A major development of the DG method was
carried out by Cockburn et al. \cite{DG2,DG3} to solve the
hyperbolic conservation laws. In the DG method, high-order accuracy
is achieved by means of high-order polynomial approximation within
each element rather than by means of wide stencils, and Runge-Kutta
method is used for the time discretization. Because only neighboring
elements interaction is included, it becomes compact and efficient
in the application on complex geometry. Recently, a correction
procedure via reconstruction framework (CPR) was developed by Wang
et al. \cite{CPR1,CPR2}. This method was originally developed to
solve one-dimensional conservation laws, under the name of flux
reconstruction \cite{CPR3,CPR4}. Under lifting collocation penalty,
the CPR framework was extended to two-dimensional triangular and
mixed grids. The CPR formulation is based on a nodal differential
form, with an element-wise continuous polynomial solution space. By
choosing certain correction functions, the CPR framework can unify
several well known methods, such as the DG, SV \cite{SD} and SD
\cite{SV} methods and lead to simplified versions of these methods,
at least for linear equations. The CPR method is compact because
only immediate face neighbors play a role in updating the solutions
in the current cell. Therefore, the complexity of implementation can
be reduced, especially for the simulation with
unstructured mesh. The main problem for the above DG-type schemes are the robustness of the method,
especially in the cases with discontinuities. It is certainly true that the use of limiters can save the DG methods
in the flow computations with discontinuities. But, the DG method is extremely sensitive to the limiters, because it is
hard to distinguish the continuous or discontinuous solution in a computation, especially with the changing of cell size.
Sometimes, the DG method can mysteriously get failure in a computation without clear reasons. Therefore, to pick up the trouble cells beforehand
becomes a practice in the DG method. After so many years' research on the DG method, it gives perfect results for the continuous flow simulations,
but seems have physical problem in its weak formulation in the discontinuous case.

 The ENO scheme was proposed by Harten et al.
\cite{ENO1, ENO2} and successfully applied to solve hyperbolic
conservation laws and other convection dominated problems. Following
the ENO scheme, WENO scheme \cite{WENO1, WENO2, WENO3} was further
developed. ENO scheme uses the smoothest stencil among several
candidates to approximate the numerical fluxes at cell interface for
high-order accuracy. At the same time, it avoids
spurious oscillations near discontinuities. Meanwhile, WENO scheme is
a convex linear combination of lower order reconstructions to obtain
a higher order approximation. WENO scheme improves upon ENO scheme
in robustness, smoothness of fluxes,  steady-state
convergence,  provable convergence properties, and more
efficiency. However, in both ENO and WENO schemes, large stencils
in the high-order reconstruction and Runge-Kutta time stepping are used, especially for the
multi-dimensional unstructured meshes \cite{WENO2}.
There are also many other high-order schemes which can be found in the literature.
The DG method for its compactness and the WENO for the reconstruction are mostly related to the current research for the development of high-order compact
gas-kinetic scheme.

The gas-kinetic scheme (GKS) based on the Bhatnagar-Gross-Krook
(BGK) model \cite{BGK-1,BGK-2,BGK-3} has been developed
systematically for the compressible flow computations \cite{GKS-Xu1,
GKS-Xu2,GKS-Kumar,GKS-Jiang,GKS-Yang}. The gas-kinetic scheme
presents a gas evolution process from a kinetic scale to a
hydrodynamic scale, where both inviscid and viscous fluxes are
recovered from moments of a single time-dependent gas distribution
function. In discontinuous shock region, the kinetic scale physics,
such as particle free transport through upwinding, takes effect to
construct a crisp and stable shock transition. The highly non-equilibrium of
the gas distribution function in the discontinuous region provides a physically consistent
mechanism for the construction of numerical shock structure.
In smooth flow region, the hydrodynamic scale physics
corresponding to the multi-dimensional central difference
discretization will contribute mainly in the kinetic flux function,
and accurate Navier-Stokes solution can be obtained once the flow
structure is well resolved. Based on the unified coordinate
transformation, a moving-mesh gas-kinetic scheme has been developed
\cite{GKS-moving1,GKS-moving2}. With the discretization of particle
velocity space, a unified gas-kinetic scheme (UGKS) has been
developed for the flow study in entire Knudsen number regimes from
rarefied to continuum ones \cite{UGKS-Xu,UGKS-Luc,Guo}. Recently, with
the incorporation of high-order initial reconstruction, high-order
gas-kinetic schemes for the inviscid and viscous flows have been
proposed in \cite{GKS-high1, GKS-high2,GKS-high0}. The flux evaluation in the
scheme is based on the time evolution of flow variables from an
initial piece-wise discontinuous polynomials around a cell
interface, where high-order spatial and temporal derivatives of a
gas distribution function are coupled nonlinearly. However, similar to
most high-order finite volume schemes, WENO-type large stencils are needed
in the reconstruction.

In this paper, a compact third-order  gas-kinetic scheme is proposed
for the compressible Euler and Navier-Stokes equations.
Different from the Riemann solver with first-order dynamics \cite{Riemann-exact,Riemann-appro},
the gas-kinetic scheme uses a time evolution solution for the flux evaluation from an initial piecewise discontinuous polynomials.
 Besides the evaluation of the time-dependent flux function across a cell interface, the high-order
gas evolution model also provides an accurate time-dependent solution of the flow variables at a cell interface.
Following the previous work for the second-order compact gas-kinetic scheme \cite{gas-compact}, the current study concentrates on the
construction of a third-order one.
The reason for the compactness of the scheme is that it not only updates the cell averaged conservative flow variables
inside each control volume, but also provides the flow variables at the cell interface at the next time level.
As a result, both cell averaged and cell interface values can be used for the high-order initial data reconstruction.
The weak formulation of DG method doesn't have such a time accurate cell interface value.
The strong solution, which follows the time evolution of flow variables
starting from a piecewise discontinuous high-order initial data, is required in the construction of the current scheme.
Due to the additional cell interface values, a compact stencil with WENO-type reconstruction can be used in the current high-order scheme.
The current scheme not only has third-order accuracy in the smooth flow regions, but also has favorable shock capturing property in the
discontinuous cases.
In comparison with other high-order schemes, the current method doesn't use the Gaussian points for the flux evaluation along the cell interface
and the multi-stage Runge-Kutta technique.

This paper is organized as follows. Section 2 introduces the BGK
equation and the multi-dimensional high-order gas-kinetic solver.
Section 3 presents the  reconstruction with compact stencil.  Section 4 presents
numerical examples to validate the current scheme. The last section
is the conclusion.

\section{BGK equation and high-order gas-kinetic solver}

\subsection{BGK equation} The
two-dimensional BGK equation can be written as
\begin{equation}\label{bgk}
f_t+uf_x+vf_y=\frac{g-f}{\tau},
\end{equation}
where $f$ is the gas distribution function, $g$ is the corresponding
equilibrium state and $\tau$ is the collision time. The collision
term satisfies the compatibility condition
\begin{equation}\label{compatibility}
\int \frac{g-f}{\tau}\psi d\Xi=0,
\end{equation}
where $\psi=(1,u,v,\displaystyle \frac{1}{2}(u^2+v^2+\xi^2))$,
$d\Xi=dudvd\xi^1...d\xi^{K}$, $K$ is the number of internal freedom,
i.e.  $K=(4-2\gamma)/(\gamma-1)$ for two-dimensional flows and
$\gamma$ is the specific heat ratio.

Based on the Chapman-Enskog expansion of the BGK model, the Euler and
Navier-Stokes, Burnett, and Super-Burnett equations can be derived
\cite{BGK-3, GKS-Xu1,Ohwada-Xu}. In the smooth region, the gas distribution
function can be expanded as
\begin{align*}
f=g-\tau D_{\textbf{u}}g+\tau D_{\textbf{u}}(\tau
D_{\textbf{u}})g-\tau D_{\textbf{u}}[\tau D_{\textbf{u}}(\tau
D_{\textbf{u}})g]+...,
\end{align*}
where $D_{\textbf{u}}=\displaystyle\frac{\partial}{\partial
t}+\textbf{u}\cdot \nabla$. By truncating different orders of
$\tau$, the corresponding macroscopic equations can be derived. For
the Euler equations, the zeroth order truncation is taken, i.e.
$f=g$. For the Navier-Stokes equations, the first order truncation
is
\begin{align}\label{ns}
f=g-\tau (ug_x+vg_y+g_t).
\end{align}
Based on the higher order truncations, the Burnett and super-Burnett
eqautions can be obtained.

Taking moments of the BGK equation Eq.\eqref{bgk} and integrating
with respect to time and space, the finite volume scheme can be
obtained
\begin{align}\label{finite}
U_{ij}^{n+1}=U_{ij}^{n}&+\frac{1}{\Delta x\Delta
y}\int_{t^n}^{t^{n+1}}\int_{-\frac{\Delta y}{2}}^{\frac{\Delta
y}{2}}[F_{i-1/2,j}(t,y)-F_{i+1/2,j}(t,y)]dy
dt\nonumber\\&+\frac{1}{\Delta x\Delta
y}\int_{t^n}^{t^{n+1}}\int_{-\frac{\Delta x}{2}}^{\frac{\Delta
x}{2}}[G_{i,j-1/2}(t,x)-G_{i,j+1/2}(t,x)]dxdt,
\end{align}
where $U=(\rho,\rho U,\rho V,\rho W,\rho E)$ are the conservative
variables, $F_{i+1/2,j}(t,y)$ and $G_{i,j+1/2}(t,x)$ are
time-dependent numerical fluxes in the $x$ and $y$ directions, which
can be obtained by taking moments of the gas distribution function
at the cell interface,
\begin{align*}
F_{i+1/2,j}(t,y)=\int u\psi f(x_{i+1/2},y,t,u,v,\xi)d\Xi.
\end{align*}
Similarly, the fluxes $G_{i,j+1/2}$ in the $y$ direction can be
obtained.

\subsection{High-order gas-kinetic solver}
To construct numerical fluxes, the integral solution of the BGK
equation Eq.\eqref{bgk} at the cell interface can be written as
\begin{align}\label{integral}
f(x_{i+1/2},y,t,u,v)=&\frac{1}{\tau}\int_0^tg(x',y',t',u,v)e^{-(t-t')/\tau}dt'+e^{-t/\tau}f_0(-ut,y-vt),
\end{align}
where $x_{i+1/2}=0$ is the location of cell interface,
$x_{i+1/2}=x'+u(t-t')$ and $y=y'+v(t-t')$  are the particle
trajectories. In the above integral solution, the initial term $f_0$
accounts for the free transport mechanism along particle
trajectories, which represents the kinetic scale physics. The
integration of equilibrium state along the particle trajectories
represents the accumulating effect of an equilibrium state, which is
related to the hydrodynamic scale flow physics. The flow behavior at
the cell interface depends on the ratio of time step and local
particle collision time $\Delta t/\tau$.

To construct a multidimensional third-order gas-kinetic solver, the
following notations are introduced firstly
\begin{align*}
a_1=&(\partial g/\partial x)/g, a_2=(\partial g/\partial y)/g,
A=(\partial g/\partial t)/g, B=(\partial A /\partial t),\\
d_{11}&=(\partial a_1/\partial x), d_{12}=(\partial a_1/\partial
y)=(\partial a_2/\partial x), d_{22}=(\partial a_2/\partial y),
\\
&b_{1}=(\partial a_1/\partial t)=(\partial A/\partial x),
b_{2}=(\partial a_2/\partial t)=(\partial A/\partial y),
\end{align*}
where $g$ is an equilibrium state. The dependence of these
coefficients on particle velocity can be expanded as the following
form \cite{GKS-Xu2}
\begin{align*}
a_1=a_{11}+a_{12}u+&a_{13}v+a_{14}\displaystyle
\frac{1}{2}(u^2+v^2+\xi^2),\\
&...\\
B=B_{1}+B_{2}u+&B_{3}v+B_{4}\displaystyle
\frac{1}{2}(u^2+v^2+\xi^2).
\end{align*}

%\begin{figure}[!h]
%\centering
%\includegraphics[width=0.45\textwidth]{0-schematic-f}
%\caption{\label{1d-schematic-f} Non-equilibrium and equilibrium
%states, and the reconstructed polynomials across the cell
%interface.}
%\end{figure}

For the kinetic part of the integral solution Eq.\eqref{integral},
the gas distribution function can be constructed as
\begin{equation}\label{f0}
f_0=f_0^l(x,y,u,v)H(x)+f_0^r(x,y,u,v)(1-H(x)),
\end{equation}
where $H(x)$ is the Heaviside function,  $f_0^l$ and $f_0^r$ are the
initial gas distribution functions on both sides of a cell interface,
which have one to one correspondence with the initially reconstructed polynomials of macroscopic
flow variables on both
sides of the cell interface. To construct a third-order scheme, the Taylor
expansion for the gas distribution function in space and time at
$(x,y)=(0,0)$ is expressed as
\begin{align*}
f_0^k(x,y)=f_G^k(0,0)&+\frac{\partial f_G^k}{\partial
x}x+\frac{\partial f_G^k}{\partial y}y+\frac{1}{2}\frac{\partial^2
f_G^k}{\partial x^2}x^2+\frac{\partial^2 f_G^k}{\partial x\partial
y}xy+\frac{1}{2}\frac{\partial^2 f_G^k}{\partial y^2}y^2,\nonumber
\end{align*}
where $k=l,r$. For the Euler equations, $f_{G}^k=g_k$ and the
kinetic part of Eq.\eqref{integral} can be obtained. For the
Navier-Stokes equations, according to Eq.\eqref{ns} and the
notations introduced above, the distribution function is
\begin{align*}
f_{G}^k=g_k-\tau(a_{1k}u+a_{2k}v+A_k)g_k,
\end{align*}
and the corresponding kinetic part of Eq.\eqref{integral} can be
written as
\begin{align}
&e^{-t/\tau}f_0^k(-ut,y-vt,u,v)\nonumber\\
=&C_7g_k[1-\tau(a_{1k}u+a_{2k}v+A_k)]\nonumber\\
+&C_8g_k[a_{1k}u-\tau((a_{1k}^2+d_{11k})u^2+(a_{1k}a_{2k}+d_{12k})uv+(A_ka_{1k}+b_{1k})u)]\nonumber\\
+&C_8g_k[a_{2k}v-\tau((a_{1k}a_{2k}+d_{12k})uv+(a_{2k}^2+d_{22k})v^2+(A_ka_{2k}+b_{2k})v)]\nonumber\\
+&C_7g_k[a_{2k}-\tau((a_{1k}a_{2k}+d_{12k})u+(a_{2k}^2+d_{22k})v+(A_ka_{2k}+b_{2k}))]y\nonumber\\
+&\frac{1}{2}C_7g_k[(a_{1k}^2+d_{11k})(-ut)^2+2(a_{1k}a_{2k}+d_{12k})(-ut)(y-vt)+(a_{2k}^2+d_{22k})(y-vt)^2],\label{dis2}
\end{align}
where $g_{k}$ are the equilibrium states at both sides of the cell
interface, and the coefficients $a_{1k},...,A_k$ are defined
according to the expansion of $g_{k}$.

After determining the kinetic part $f_0$, the equilibrium state $g$
in the integral solution Eq.\eqref{integral} can be constructed
consistent with $f_0$ as follows
\begin{align}\label{equli}
g=g_0+\frac{\partial g_0}{\partial x}x+&\frac{\partial g_0}{\partial
y}y+\frac{\partial g_0}{\partial t}t+\frac{1}{2}\frac{\partial^2
g_0}{\partial x^2}x^2+\frac{\partial^2 g_0}{\partial x\partial
y}xy+\frac{1}{2}\frac{\partial^2 g_0}{\partial
y^2}y^2\nonumber\\
&+\frac{1}{2}\frac{\partial^2 g_0}{\partial t^2}t^2+\frac{\partial^2
g_0}{\partial x\partial t}xt+\frac{\partial^2 g_0}{\partial
y\partial t}yt,
\end{align}
where $g_{0}$ is the equilibrium state located at interface, which
can be determined through the compatibility condition
Eq.\eqref{compatibility}
\begin{align}\label{compatibility2}
\int\psi g_{0}d\Xi=U_0=\int_{u>0}\psi g_{l}d\Xi+\int_{u<0}\psi
g_{r}d\Xi.
\end{align}
Based on Taylor expansion for the equilibrium state
Eq.\eqref{equli}, the hydrodynamic part in Eq.\eqref{integral} can
be written as
\begin{align}\label{dis1}
\frac{1}{\tau}\int_0^t
g&(x',y',t',u,v)e^{-(t-t')/\tau}dt'\nonumber\\
=&C_1g_0+C_2g_0\overline{a}_1u+C_2g_0\overline{a}_2v+C_1g_0\overline{a}_2y+C_3g_0\overline{A}\nonumber\\
+&\frac{1}{2}C_4g_0(\overline{a}_1^2+\overline{d}_{11})u^2+C_6g_0(\overline{A}\overline{a}_1+\overline{b}_{1})u+\frac{1}{2}C_5g_0(\overline{A}^2+\overline{B})\nonumber\\
+&\frac{1}{2}C_1g_0(\overline{a}_2^2+\overline{d}_{22})y^2+C_2g_0(\overline{a}_2^2+\overline{d}_{22})vy+\frac{1}{2}C_4g_0(\overline{a}_2^2+\overline{d}_{22})v^2\nonumber \\
+&C_2g_0(\overline{a}_1\overline{a}_2+\overline{d}_{12})uy+C_4g_0(\overline{a}_1\overline{a}_2+\overline{d}_{12})uv\nonumber\\
+&C_3g_0(\overline{A}\overline{a}_2+\overline{b}_{2})y+C_6g_0(\overline{A}\overline{a}_2+\overline{b}_{2})v,
\end{align}
where the coefficients
$\overline{a}_1,\overline{a}_2,...,\overline{A},\overline{B}$ are
defined from the expansion of the equilibrium state $g_0$.
The coefficients $C_i, i=1,...,8$ in Eq.\eqref{dis1} and Eq.\eqref{dis2}
are given by
\begin{align*}
C_1=1-&e^{-t/\tau}, C_2=(t+\tau)e^{-t/\tau}-\tau, C_3=t-\tau+\tau e^{-t/\tau},C_4=-(t^2+2t\tau)e^{-t/\tau},\\
&C_5=t^2-2t\tau,C_6=-t\tau(1+e^{-t/\tau}),C_7=e^{-t/\tau},C_8=-te^{-t/\tau}.
\end{align*}
Substituting Eq.\eqref{dis1} and Eq.\eqref{dis2} into the integral
solution Eq.\eqref{integral}, the gas distribution function at the
cell interface can be obtained.

For the smooth flow, the polynomials at both sides of the cell
interface take the same polynomial $\overline{U}(x)$, which gives
$g_0=g_l=g_r$ and identical slopes. Consequently, the gas
distribution function Eq.\eqref{integral} will reduce to the
continuous one
\begin{align}\label{dis3}
f=&g_0[1-\tau(\overline{a}_1u+\overline{a}_2v+\overline{A})]\nonumber\\
+&g_0[\overline{a}_2-\tau((\overline{a}_1\overline{a}_2+\overline{d}_{12})u+(\overline{a}_2^2+\overline{d}_{22})v+
(\overline{A}\overline{a}_2+\overline{b}_2))]y\nonumber\\
+&g_0[\overline{A}-\tau((\overline{A}\overline{a}_1+\overline{b}_1)
u+(\overline{A}\overline{a}_2+\overline{b}_2)v+(\overline{A}^2+\overline{B}))]t\nonumber\\
+&g_0[\frac{1}{2}(\overline{a}_2^2+\overline{d}_{22})y^2+(\overline{A}\overline{a}_2+\overline{b}_2)yt+\frac{1}{2}(\overline{A}^2+\overline{B})t^2].
\end{align}

The superscripts or subscripts of the coefficients $a_1, a_2,...,A,
B$ in Eq.\eqref{dis2}, Eq.\eqref{dis1} and Eq.\eqref{dis3} are
omitted for simplicity and they are determined by the spatial
derivatives of macroscopic flow variables and the compatibility
condition \cite{GKS-high2} as follows
\begin{align}
\begin{cases}
&\displaystyle\langle a_1\rangle =\frac{\partial U}{\partial x},
\langle a_2\rangle =\frac{\partial U }{\partial y},  \langle
A+a_1u+a_2v \rangle=0,\\ &\displaystyle\langle a_1
^2+d_{11}\rangle=\frac{\partial^2 U }{\partial x^2}, \langle a_2
^2+d_{22}\rangle=\frac{\partial^2 U }{\partial y^2}, \langle
a_1a_2+d_{12}\rangle=\frac{\partial^2
U}{\partial x\partial y},\\
&\displaystyle\langle(a_1 ^2+d_{11})u+(a_1a_2+d_{12})v+(Aa_1+b_1)\rangle=0,\\
&\displaystyle\langle(a_1a_2+d_{12})u+(a_2 ^2+d_{22})v+(Aa_2+b_2)\rangle=0,\\
&\displaystyle\langle(Aa_1+b_1)u+(Aa_2+b_2)v+(A^2+B)\rangle=0,
\end{cases}
\end{align}
where $<...>$ are the moments of gas distribution function, and
defined by
\begin{align*}
<...>=\int g(...)\psi d\Xi.
\end{align*}

\section{Compact reconstruction}

In the traditional high-order schemes, a high-order polynomial is
reconstructed or updated inside each cell, and the exact Riemann solver
\cite{Riemann-exact} or approximate Riemann solvers
\cite{Riemann-appro} are used to provide flux function at the cell
interface. The Riemann solver presents the Euler solution from two constant states, which has the wave
propagation in the normal direction of the cell interface only.
Due to its constant flux and state at the cell interface, the cell interface solution has only first-order accuracy, which can be hardly used
in the reconstruction at the beginning of next time step.
In the gas-kinetic scheme,
besides the numerical fluxes, the pointwise values of the
macroscopic variables at a cell interface can be obtained as well by taking moments of the gas
distribution function \cite{gas-compact},
\begin{align}\label{point}
U_{i+1/2,j}(t,y)&=\int \psi f(x_{i+1/2},y,t,u,v,\xi)d\Xi.
\end{align}
As shown in the last section, the whole curve of the polynomial of
the macroscopic variables will participate the flow evolution, and the spatial and
temporal derivatives of the gas distribution function are coupled
nonlinearly. This pointwise values at the cell interface
(Eq.\eqref{point}) is a strong high-order dynamic solution, which can be used in the reconstruction stage at the
beginning of next time step.
This is also the main point we would like to emphasize in this paper that the low order dynamics of the Riemann solver may be the bottle
neck for the development of high-order compact numerical schemes. The use of the weak solution, such as DG, is mainly to avoid the use of high-order flow dynamics.
In the following subsections, a third-order compact reconstruction will
be presented for both one and two dimensional cases, in which the
pointwise values at the cell interface and the cell averaged values in the neighboring cells only are used in the high-order reconstruction.

\begin{figure}[!h]
\centering
\includegraphics[width=0.75\textwidth]{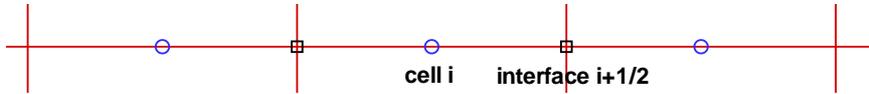}
\caption{\label{1d-schematic} One dimensional stencil for the cell
$I_{i}$. The circles represent the cell averaged values
$\overline{U}_{i-1},\overline{U}_i,\overline{U}_{i+1}$ and the
squares represent the pointwise values $U_{i-1/2}, U_{i+1/2}$ at the
cell interface.}
\end{figure}

\subsection{One-dimensional initial data reconstruction}
In the one-dimensional case, the stencils of the compact
reconstruction for cell $I_i$ are shown in Fig.\ref{1d-schematic}.
Two types of stencils are considered in the reconstruction and the
procedure are listed as follows
\begin{enumerate}
\item For one-sided stencils
$S_1=\{\overline{U}_{i-1},U_{i-1/2},\overline{U}_i\}$ and
$S_2=\{\overline{U}_{i},U_{i+1/2},\overline{U}_{i+1}\}$, a
quadratic polynomial can be defined as
$\phi_k(x)=\overline{U}_{i}+a_2^kx+a_3^k\omega_x, k=1,2$, where
$\displaystyle\omega_x=\frac{1}{2}(x^2-\frac{1}{12}\Delta x^2)$.
According to the cell averaged values $\overline{U}_{i+(-1)^{k}}$
and point values $U_{i+(-1)^{k}/2}$, the coefficients in $\phi_1,
\phi_2$ can be expressed as
\begin{align*}
&a_{3}^k=3(\overline{U}_{i}+\overline{U}_{i+(-1)^{k}}-2U_{i+(-1)^{k}/2}))/\Delta x^2,\\
a_{2}^k=(-1)^{k}[2(U&_{i+(-1)^{k}/2}-\overline{U}_{i})/\Delta
x-a_{3}^k\Delta x/6]=(-1)^{k}(\widetilde{a}_{2}^k-a_{3}^k\Delta
x/6).
\end{align*}
To deal with possible flow with discontinuity, the above coefficients are
limited as
\begin{align}\label{limiter}
\left\{\begin{aligned}
         &a_{3}=minmod\{a_{3}^1,a_{3}^2\},\\
         &a_{2}=minmod\{-\widetilde{a}_{2}^1+a_{3}\Delta
x/6,\widetilde{a}_{2}^2-a_{3}\Delta x/6\},
                          \end{aligned} \right.
\end{align}
where $minmod\{\cdot,\cdot\}$ is the minmod limiter. With these
modified coefficients, the polynomial $\displaystyle
U_1(x)=\overline{U}_{i}+a_2x+a_3\omega_x$ from two
one-sided stencils can be fully determined.
\item For the central stencil $\displaystyle
S_3=\{\overline{U}_{i-1},\overline{U}_i,\overline{U}_{i+1}\}$, the
polynomial $\displaystyle U_2(x)$ can be obtained according to
\begin{align*}
\int_{I_{i-1}}U_2(x)dx=\overline{U}_{i-1},
\int_{I_{i+1}}U_2(x)dx=\overline{U}_{i+1}.
\end{align*}
\item With the polynomials $U_i(x), i=1,2$ corresponding to the
one-sided and central stencils, the non-linear weights $\omega_i$
\cite{WENO2} are used to construct the combined polynomial in the cell
$I_{i}$ as follows
\begin{align}\label{average}
U(x)=\omega_1U_1(x)+\omega_2U_2(x),
\end{align}
where $\displaystyle\omega_i=\frac{\alpha_i}{\sum_i\alpha_i}$,
$\displaystyle\alpha_i=\frac{1}{(\varepsilon+I(U_i))^2}$,
$\varepsilon$ is a small number and $I(U_i)$ is the smooth indicator
of $U_i$, which is expressed as
\begin{align*}
I(U_i)=\sum_{1\leq l\leq 2}\int_{I_i}h^{2l-1}(U_i^{(l)})^2dx,
\end{align*}
where $U_i^{(l)}$ is the $l$-th order derivative of $U_i$.\\
In most cases, the polynomial $U_1(x)$ corresponding to the
one-sided stencils can well resolve the discontinuities, and the
central stencil is introduced to improve the accuracy in the smooth
region. In the computation, the smooth indicator of $U_1$ directly
takes $I(U_1)=Ch^2$. Thus, $I(U_2)=O(1)$ in the region with
discontinuity and $U_1$ is the dominant one in Eq.\eqref{average};
in the smooth region, $I(U_2)=O(h^2)$ and $U(x)$ is the average of
$U_1(x)$ and $U_2(x)$.
\end{enumerate}

\begin{figure}[!h]
\centering
\includegraphics[width=0.43\textwidth]{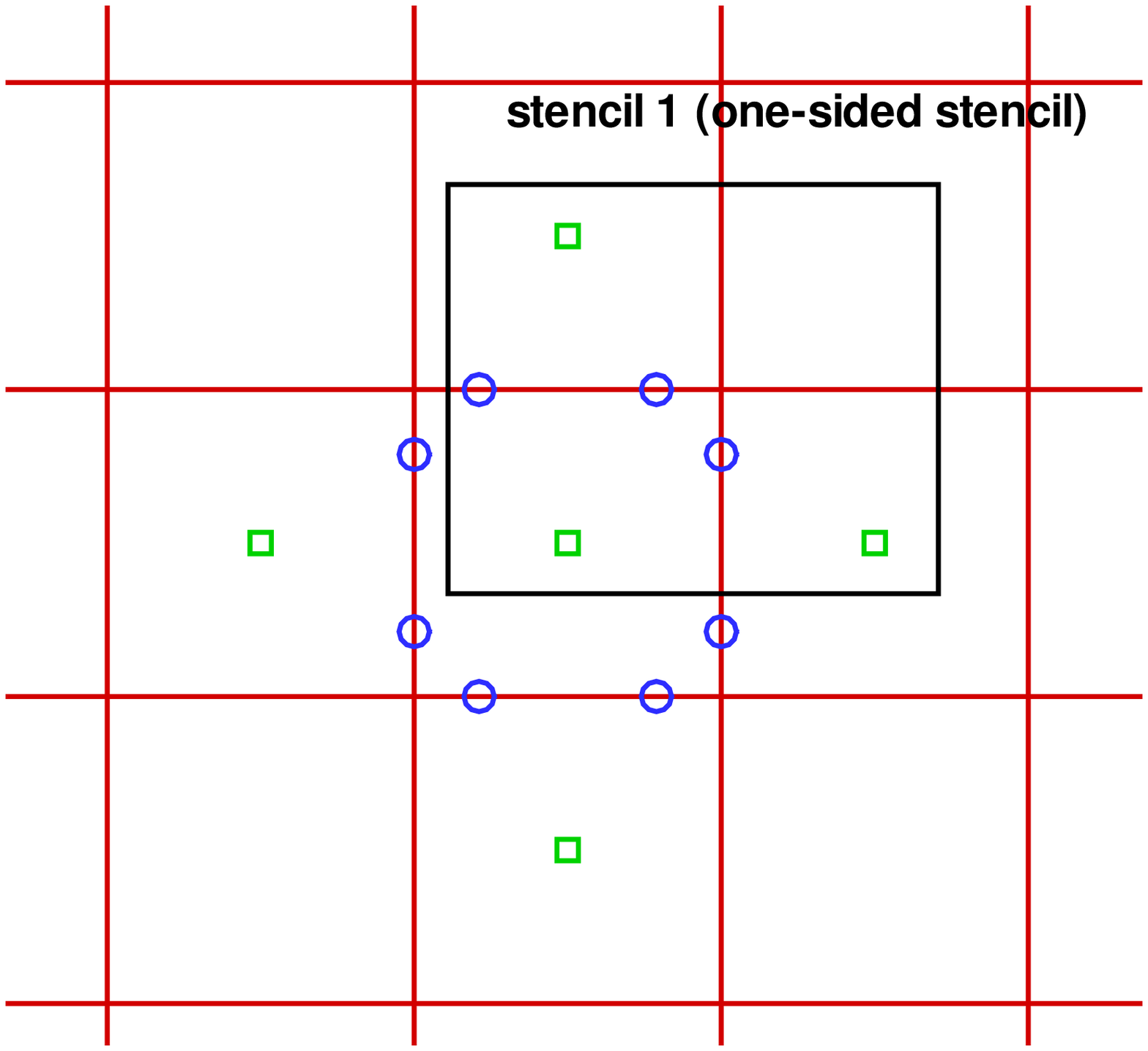}
\includegraphics[width=0.43\textwidth]{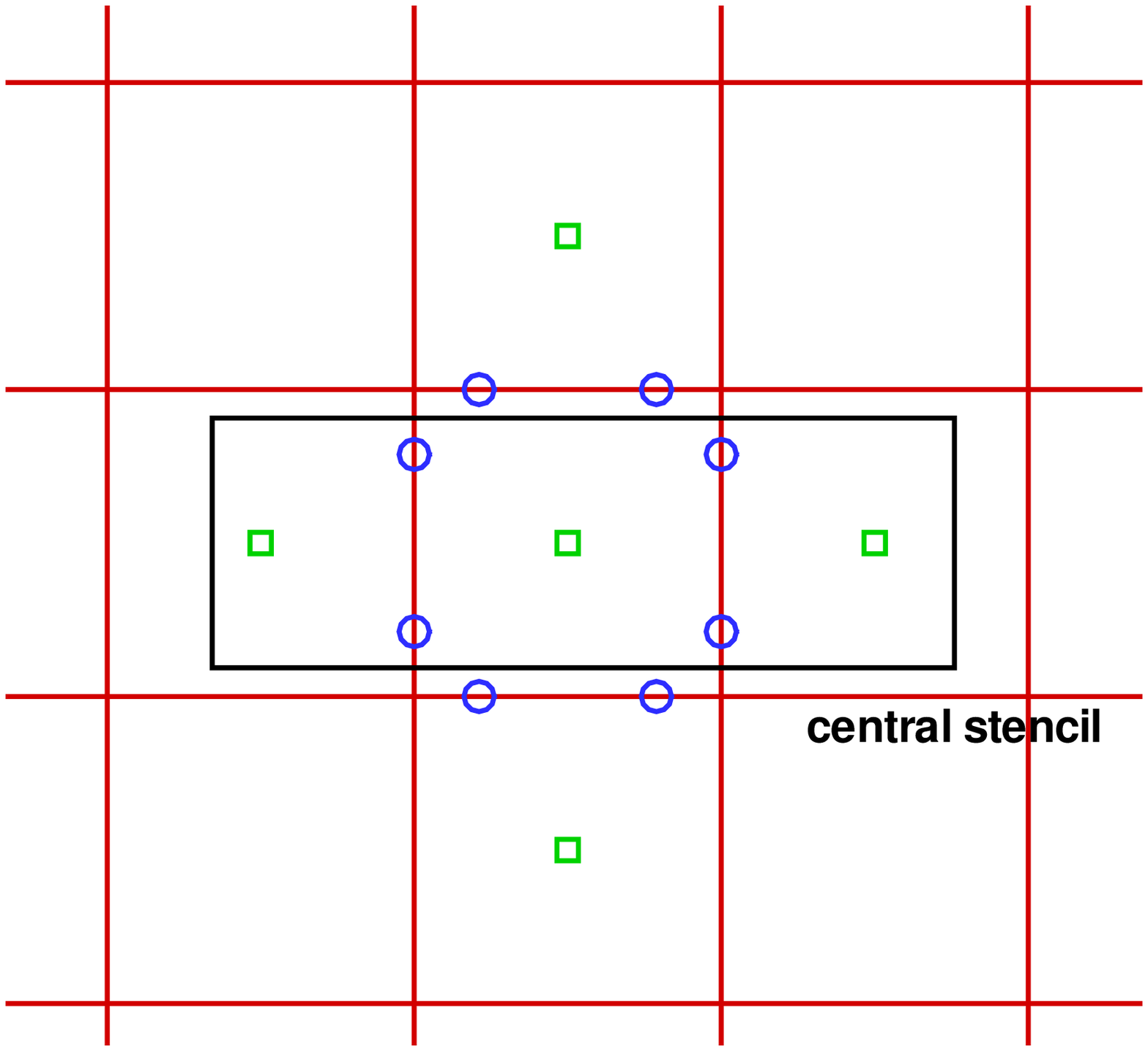}
\caption{\label{schematic-2d} Two types of stencils for the compact
reconstruction of cell $I_{ij}$. Left one has four one-sided stencils,
where only stencil 1 is shown; the right one is the central stencil.
The squares represent the cell averaged values of each cell; circles
are the pointwise values at the Gaussian points, which can be obtained from the solution in Eq.(\eqref{point}).}
\end{figure}

\subsection{Two-dimensional initial data reconstruction}
Similar to the one-dimensional case, two types of stencils are used in
the data reconstruction i.e. one-sided stencils and central stencil. For
the cell $I_{ij}$, the quadratic polynomial $\phi_k(x,y)$ are
defined by
\begin{align*}
\phi_k(x,y)=\overline{U}_{ij}+a_{2}^kx+a_{3}^ky+a_4^k\omega_x+a_5^k\omega_y+a_6^kxy,
\end{align*}
where $\overline{U}_{ij}$ is the cell averaged value of $I_{ij}$.
For the one-sided stencil 1, the polynomial $\phi_1(x,y)$ can be
determined according to three cell averaged values and four pointwise
values at Gaussian integration points as follows
\begin{align*}
\int_{I_{i+1,j}}\phi_1(x,y)dxdy=\overline{U}_{i+1,j},&~~
\int_{I_{i,j+1}}\phi_1(x,y)dxdy=\overline{U}_{i,j+1},\\
\phi_1(x_{i+\frac{1}{2}},y_{j+\frac{\sqrt{3}}{6}})=U_{i+\frac{1}{2},j^+},&~~
\phi_1(x_{i+\frac{1}{2}},y_{j-\frac{\sqrt{3}}{6}})=U_{i+\frac{1}{2},j^-},\\
\phi_1(x_{i+\frac{\sqrt{3}}{6}},y_{j+\frac{1}{2}})=U_{i^+,j+\frac{1}{2}},&~~
\phi_1(x_{i-\frac{\sqrt{3}}{6}},y_{j+\frac{1}{2}})=U_{i^-,j+\frac{1}{2}}.
\end{align*}
The least square solution for the above over-determined system are
written as
\begin{align}\label{over-sys}
\begin{cases}
a_{4}^1&=3[(\overline{U}_{i,j}+\overline{U}_{i+1,j})-(U_{i+\frac{1}{2},j^+}+U_{i+\frac{1}{2},j^-})]/\Delta x^2,\\
a_{2}^1&=(U_{i+\frac{1}{2},j^+}+U_{i+\frac{1}{2},j^-}-2\overline{U}_{ij})/\Delta
x-a_{4}^1\Delta x/6,\\
a_{5}^1&=3[(\overline{U}_{i,j}+\overline{U}_{i,j+1})-(U_{i^+,j+\frac{1}{2}}+U_{i^-,j+\frac{1}{2}})]/\Delta y^2\\
a_{3}^1&=(U_{i^+,j+\frac{1}{2}}+U_{i^-,j+\frac{1}{2}}-2\overline{U}_{ij})/\Delta
y-a_{5}^1\Delta y/6,\\
a_{6}^1&=\sqrt{3}[U_{i+\frac{1}{2},j^+}-U_{i+\frac{1}{2},j^-}+U_{i^+,j+\frac{1}{2}}-U_{i^-,j+\frac{1}{2}}]/\Delta
x\Delta y-(a_{2}^1/\Delta y+a_{3}^1/\Delta x).
\end{cases}
\end{align}
Similarly, the coefficients $a^k_m$ of the quadratic polynomials
$\phi_k(x,y), k=2,3,4, m=2,...,6$ corresponding to other three
one-sided stencils can be also obtained.

Due to the decoupling of the normal and tangential derivatives
(Eq.\eqref{over-sys}) for the quadratic polynomial $\phi_k(x,y)$, the
normal derivatives $a_2, a_4$ and tangential derivatives $a_3, a_5$
can be modified according to the limiter Eq.\eqref{limiter},
respectively. With the modified coefficients $a_2, a_3$, the cross
derivative $a_6^1$ is modified as
\begin{align*}
\widetilde{a}_{6}^1=\sqrt{3}[U_{i+\frac{1}{2},j^+}-U_{i+\frac{1}{2},j^-}+U_{i^+,j+\frac{1}{2}}-U_{i^-,j+\frac{1}{2}}]/\Delta
x\Delta y-(a_{2}/\Delta y+a_{3}/\Delta x).
\end{align*}
Choosing the one with the smallest absolute value from
$\widetilde{a}_6^k$, $k=1,...,4$, the limiting procedure is done and
a quadratic polynomial
$U_1(x,y)=\overline{U}_{i,j}+a_{2}x+a_{3}y+a_4\omega_x+a_5\omega_y+a_6xy$
can be obtained for the four one-sided stencils.

For the central stencil $S_5$, the polynomial $\displaystyle
U_2(x,y)$ can be obtained as follows
\begin{align*}
\int_{I_{i-1}}U_2(x)dx=\overline{U}_{i-1},&~~
\int_{I_{i+1}}U_2(x)dx=\overline{U}_{i+1},\\
U_2(x_{i+\frac{1}{2}},y_{j+\frac{\sqrt{3}}{6}})=U_{i+\frac{1}{2},j^+},&~~
U_2(x_{i+\frac{1}{2}},y_{j-\frac{\sqrt{3}}{6}})=U_{i+\frac{1}{2},j^-},\\
U_2(x_{i-\frac{1}{2}},y_{j+\frac{\sqrt{3}}{6}})=U_{i-\frac{1}{2},j^+},&~~
U_2(x_{i-\frac{1}{2}},y_{j-\frac{\sqrt{3}}{6}})=U_{i-\frac{1}{2},j^-}.
\end{align*}

With the non-linear weights for the two-dimensional reconstruction,
the polynomial in the cell $I_{ij}$ is constructed as
\begin{align*}
U(x)=\omega_1U_1(x)+\omega_2U_2(x),
\end{align*}
where $\omega_i, i=1,2$ are the non-linear weights and details can
be found in \cite{WENO3}.

\begin{figure}[!h]
\centering
\includegraphics[width=0.35\textwidth]{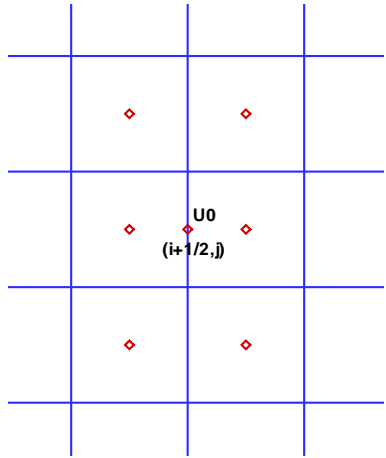}
\caption{\label{schematic2} The stencil for the equilibrium
distribution around the cell interface $(i+1/2,j)$.}
\end{figure}

\subsection{Reconstruction for equilibrium part}
In this subsection, a quadratic polynomial $\overline{U}(x,y)$
corresponding to equilibrium (hydrodynamic) part will be constructed, which is
expressed as
\begin{align*}
\overline{U}(x,y)=U_{0}+b_{2}x+b_{3}y+\frac{1}{2}b_4x^2+\frac{1}{2}b_5y^2+b_6xy.
\end{align*}
The conservative variables $U_0$ at the center of cell interface are
obtained according to the compatibility condition
Eq.\eqref{compatibility2}, in which $g_l, g_r$ are equilibrium
states corresponding to the initially reconstructed conservative variables
$U_{i+1/2}^l, U_{i+1/2}^r$ at both sides of cell interface.

To determine this polynomial with the compact stencils, six cell
averaged values are used as shown in Fig.\ref{schematic2} with the
following conditions
\begin{align*}
\iint_{I_{i+m,j+n}}\overline{U}(x,y)dxdy=\overline{U}_{i+m,j+n},
\end{align*}
where $m=0,1, n=-1,0,1$. The coefficients of  $\overline{U}(x,y)$
can be obtained by the least square procedure, and they are
expressed as
\begin{align*}
&b_2=((\overline{U}_{i+1,j+1}-\overline{U}_{i,j+1})+2(\overline{U}_{i+1,j}-\overline{U}_{ij})+(\overline{U}_{i+1,j-1}-\overline{U}_{i,j-1}))/4\Delta x,\\
&b_3=(\overline{U}_{i+1,j+1}-\overline{U}_{i+1,j-1}+\overline{U}_{i,j+1}-\overline{U}_{i,j-1})/4\Delta y,\\
&b_4=(26(\overline{U}_{i+1,j}+\overline{U}_{ij})-(\overline{U}_{i+1,j+1}+\overline{U}_{i,j+1}+\overline{U}_{i+1,j-1}+\overline{U}_{i,j-1})-48U_0)/8\Delta x^2,\\
&b_5=(\overline{U}_{i+1,j+1}-2\overline{U}_{i+1,j}+\overline{U}_{i+1,j-1}+\overline{U}_{i,j+1}-2\overline{U}_{i,j}+\overline{U}_{i,j-1})/2\Delta y^2,\\
&b_6=(\overline{U}_{i+1,j+1}-\overline{U}_{i,j+1}-\overline{U}_{i+1,j-1}+\overline{U}_{i,j-1})/2\Delta
x\Delta y.
\end{align*}

Similarly,  the quadratic polynomial
$\overline{U}(x)$ across the cell interface can be also constructed in the one-dimensional case.

\noindent{\bf{Remark:}}
In this section, a compact simple reconstruction is presented, which can be easily extended to unstructured
mesh. Theoretically, for the third-order gas-kinetic scheme, three independent
pointwise values at a cell interface (Eq.\eqref{point}) can be obtained and used for
the spatial data reconstruction. There are many choices for
the reconstruction. To obtain optimal and robust reconstruction scheme specifically to the kinetic formulation is
an interesting open question and needs further investigation.

\section{Numerical tests}

In this section, numerical tests for both inviscid flow and viscous
flow will be presented to validate our numerical scheme. For the
inviscid flow, the collision time $\tau$ takes
\begin{align*}
\tau=\epsilon \Delta t+C\displaystyle|\frac{p_l-p_r}{p_l+p_r}|\Delta
t,
\end{align*}
where $\varepsilon=0.05$ and $C=1$. For the viscous flow, we have
\begin{align*}
\tau=\frac{\mu}{p}+\displaystyle|\frac{p_l-p_r}{p_l+p_r}|\Delta t,
\end{align*}
where $p_l$ and $p_r$ denotes the pressure on the left and right sides
of the cell interface, $\mu$ is the viscous coefficient and $p$ is
the pressure at the cell interface. In the smooth flow regions, it will
reduce to $\tau=\mu/p$.
For diatomic molecules with $\gamma = 1.4$, the current gas-kinetic scheme solves the NS equations with the inclusion of bulk viscosity \cite{xu-liu-jiang}.
For monatomic gas with $\gamma = 5/3$, there is no bulk viscosity involved.
$\Delta t$ is the time step which is
determined according to the CFL condition. In the numerical tests,
CFL number takes 0.2.

For the smooth flow, the compact reconstruction is based on the
conservative variables directly; for the flow with discontinuity,
this reconstruction is based on the characteristic variables.

\subsection{Accuracy tests}
We consider two test cases to verify the numerical order of the
compact gas-kinetic scheme for the invicid flow. The first case is
the advection of density perturbation, and the initial condition is
set as follows
\begin{align*}
\rho(x)=1+0.2\sin(\pi x), U(x)=1, p(x)=1, x\in[0,2].
\end{align*}
The periodic boundary condition is adopted and thus the analytic
solution is
\begin{align*}
\rho(x)=1+0.2\sin(\pi(x-t)), U(x)=1, p(x)=1.
\end{align*}
In the computation, the uniform mesh is used, and the $L^1$ and
$L^2$ errors and orders at $t=2$ are presented in Table.\ref{tab1},
which shows the third-order accuracy.

\begin{table}[!h]
\begin{center}
\def\temptablewidth{0.8\textwidth}
{\rule{\temptablewidth}{0.5pt}}
\begin{tabular*}{\temptablewidth}{@{\extracolsep{\fill}}c|cc|cc}
mesh & $L^1$ norm & order ~ & $L^2$ norm & order  \\
\hline
50  & 6.994400E-006 &  ~~         & 7.765728E-006&    ~~\\
100 & 7.925999E-007 &  3.141535   & 8.803306E-007&    3.141004\\
200 & 1.069000E-007 &  2.890331   & 1.182243E-007&    2.896518\\
400 & 1.329999E-008 &  3.006764   & 1.516575E-008&    2.962638
\end{tabular*}
{\rule{\temptablewidth}{0.5pt}}
\end{center}
\vspace{-4mm} \caption{\label{tab1} Space accuracy test for the
advection of density perturbation.}
\end{table}

The second one is isotropic vortex propagation problem
\cite{Case-Shu}. The mean flow is $(\rho,u,v,p)=(1,1,1,1)$, and an
isotropic vortex is added to the mean flow, i.e., with perturbation
in $u, v$ and temperature $T=p/\rho$, and no perturbation in entropy
$S=p/\rho^\gamma$. The perturbation is given by
\begin{align*}
&(\delta u,\delta v)=\frac{\epsilon}{2\pi}e^{\frac{(1-r^2)}{2}}(-y,x),\\
\delta& T=-\frac{(\gamma-1)\epsilon^2}{8\gamma\pi^2}e^{1-r^2},\delta
S=0,
\end{align*}
where $r^2=x^2+y^2$ and the vortex strength $\epsilon=5$. The
computational domain is $[-5,5]\times[-5,5]$ and the periodic
boundary conditions are imposed on the boundaries in both $x$ and
$y$ directions. The exact solution is the perturbation which
propagates with the velocity $(1,1)$.  The $L^1$ and $L^2$ errors
and orders after one time period with $t=10$ are presented in
Table.\ref{tab2}, which shows that the third-order accuracy can be also
achieved.

\begin{figure}[!h]
\centering
\includegraphics[width=0.39\textwidth]{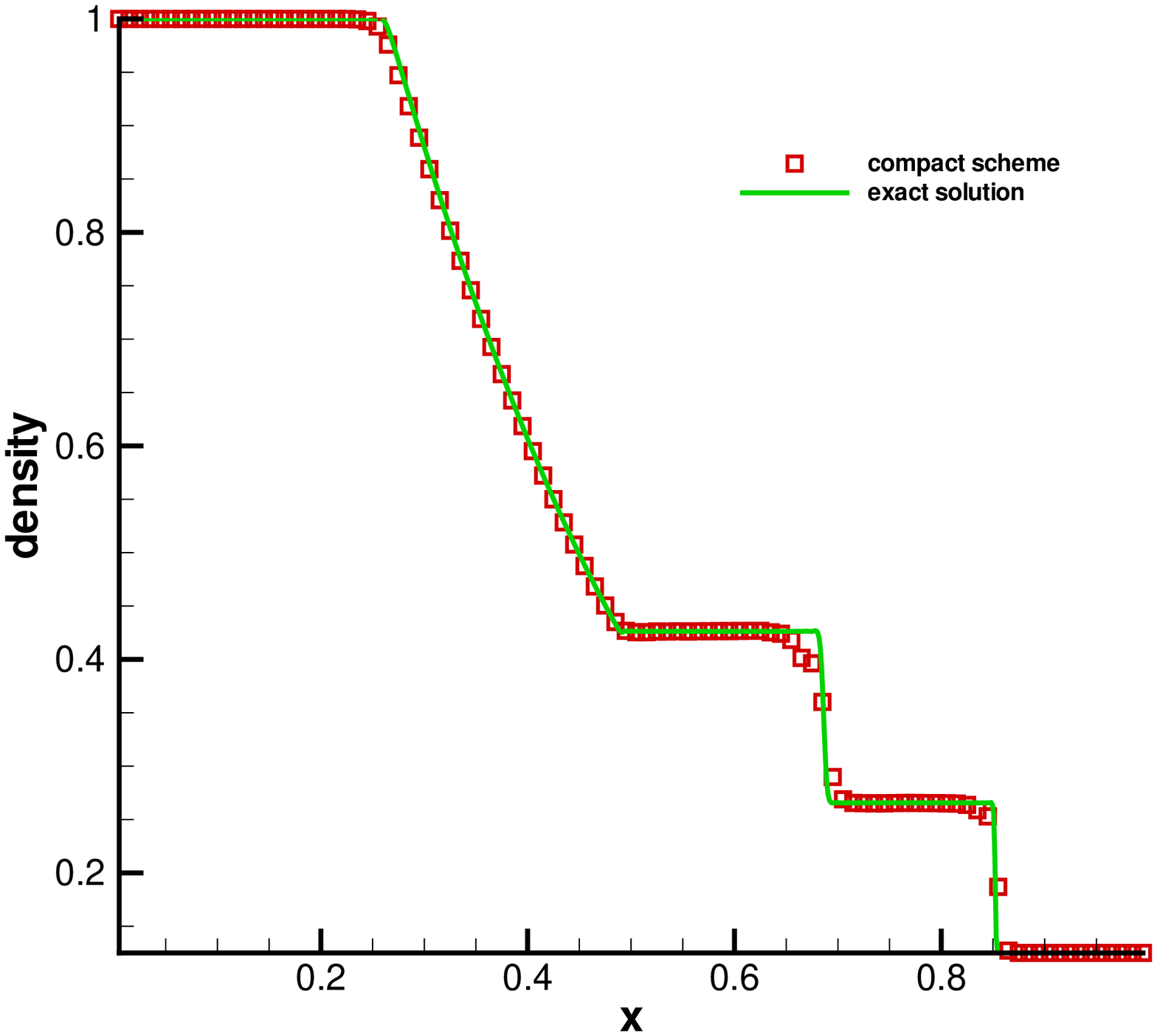}\includegraphics[width=0.39\textwidth]{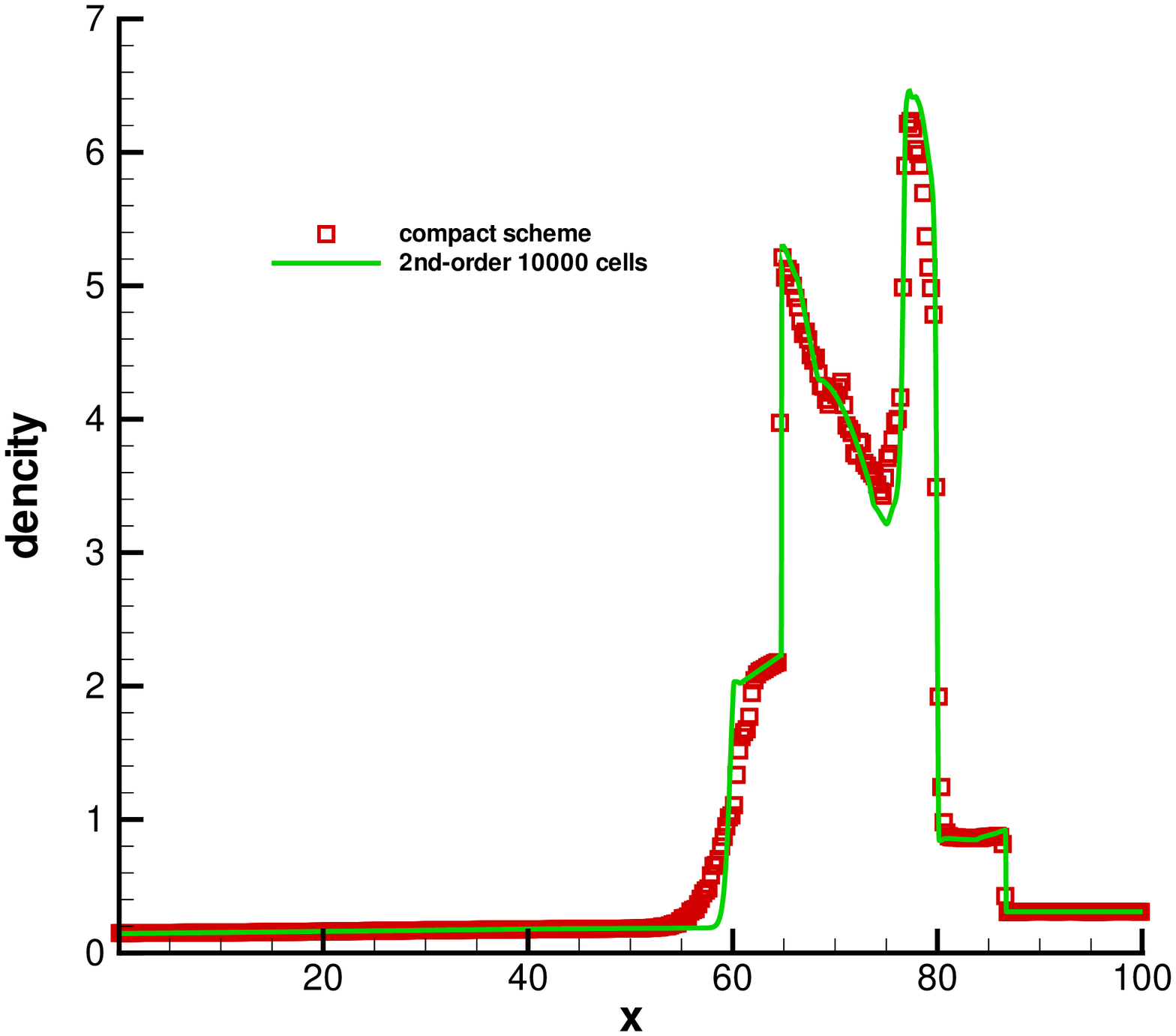}
\includegraphics[width=0.39\textwidth]{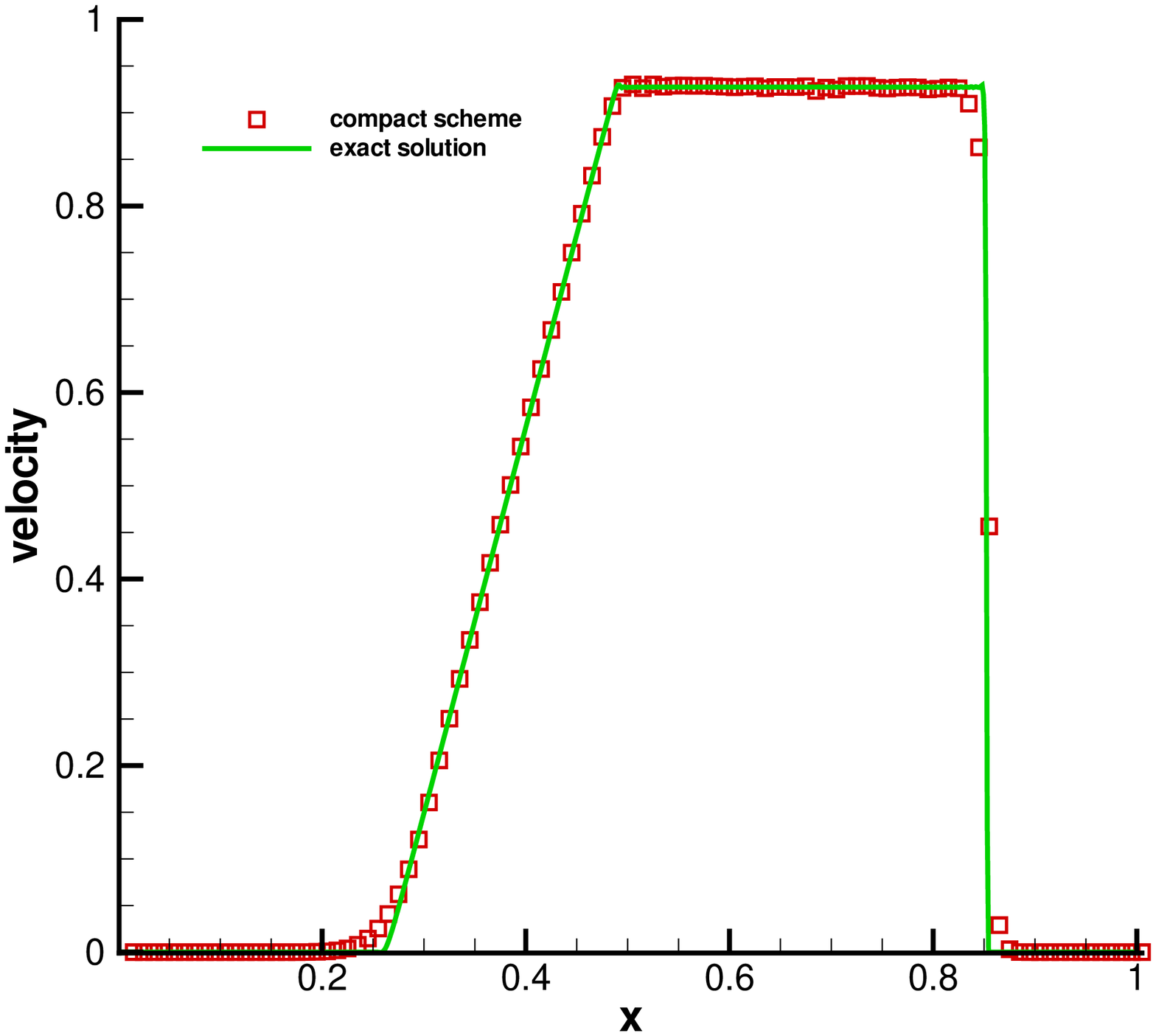}\includegraphics[width=0.39\textwidth]{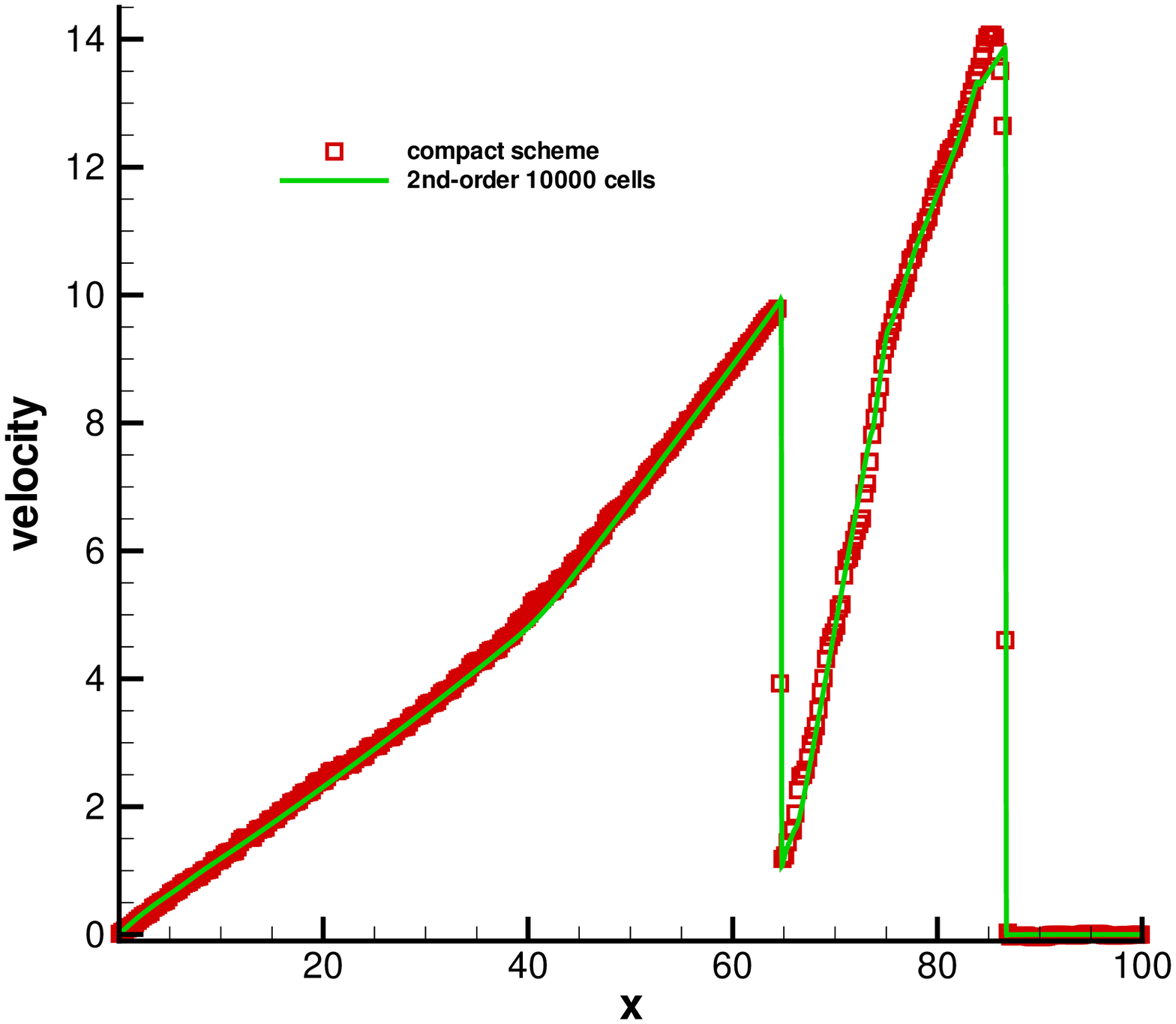}
\includegraphics[width=0.39\textwidth]{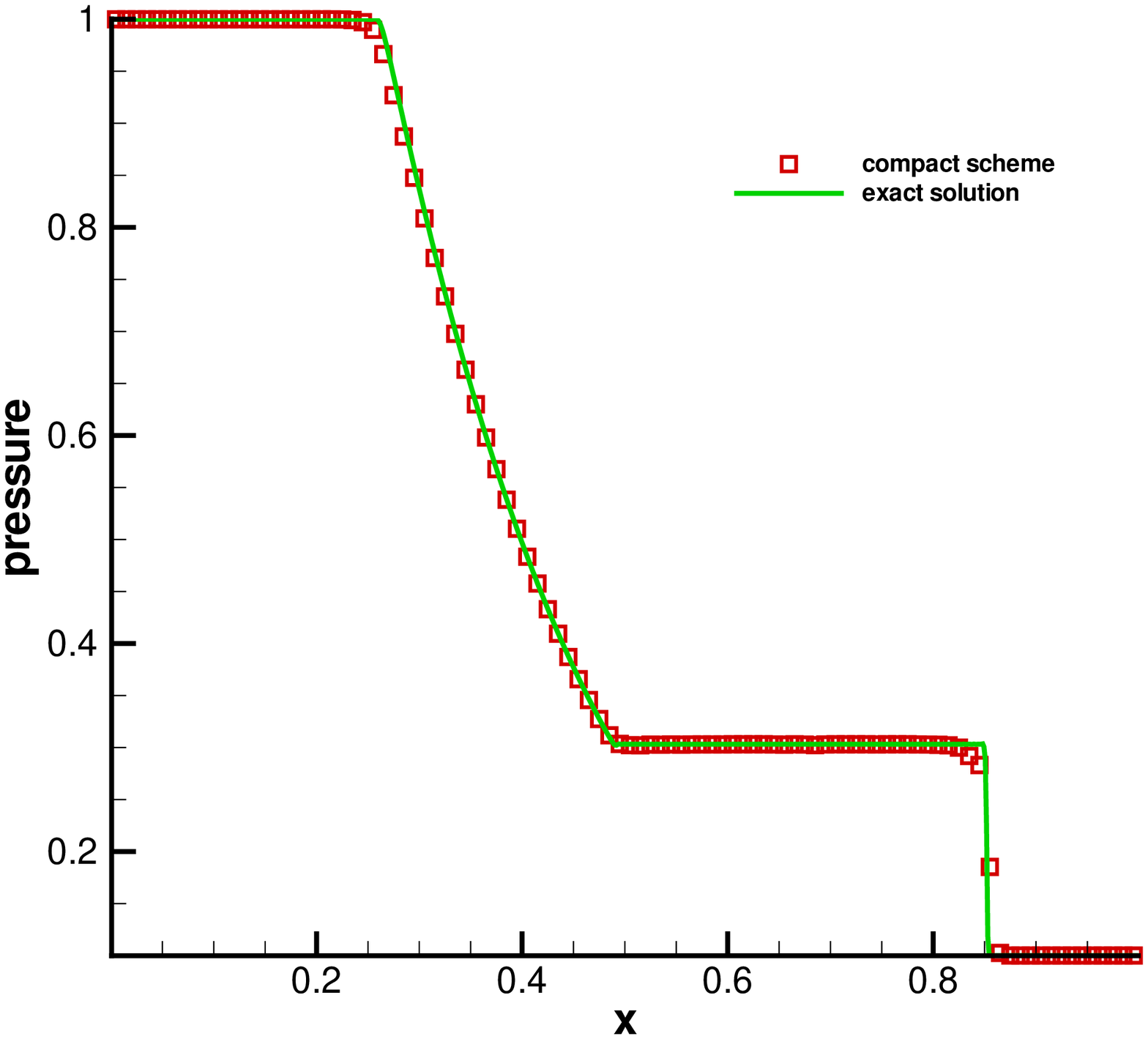}\includegraphics[width=0.39\textwidth]{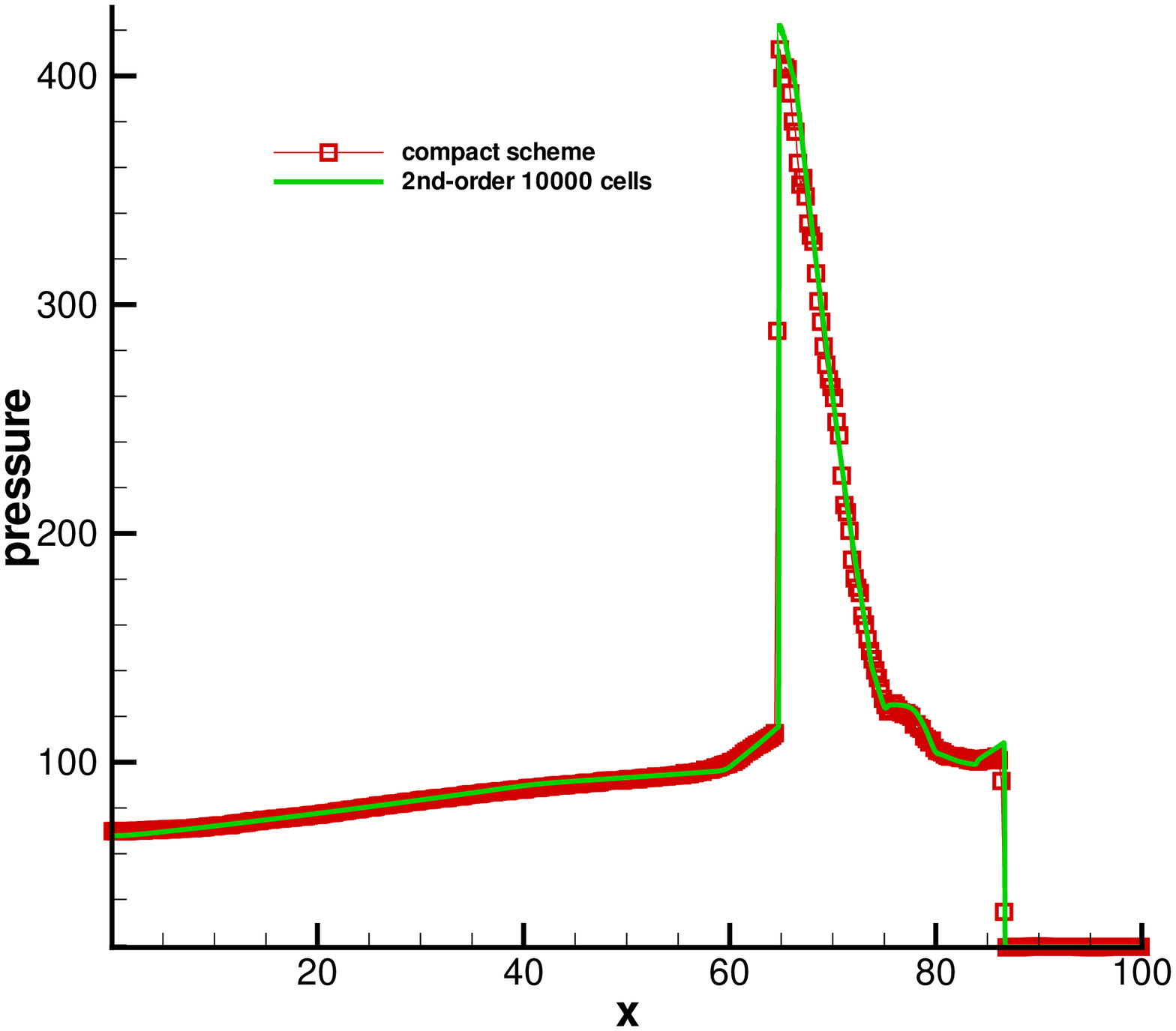}
\caption{\label{1d} Sod problem (left): the density, velocity and
pressure distributions at $t=0.2$. Blast wave problem (right): the
density, velocity and pressure distributions at $t=3.8$.}
\end{figure}

\begin{table}[!h]
\begin{center}
\def\temptablewidth{0.8\textwidth}
{\rule{\temptablewidth}{0.5pt}}
\begin{tabular*}{\temptablewidth}{@{\extracolsep{\fill}}c|cc|cc}
mesh & $L^1$ norm & order ~ & $L^2$ norm & order  \\
\hline
$21\times21$   & 3.2800781E-02   &  ~        & 6.2543898E-02 & ~\\
$41\times41$   & 5.5260700E-03   & 2.661948  & 9.4291484E-03 &2.827986 \\
$81\times81$   & 8.7312283E-04   & 2.709969  & 1.5213625E-03 &2.679191\\
$121\times121$ & 2.6716280E-04   & 2.950650  & 4.8235722E-04
&2.862094
\end{tabular*}
{\rule{\temptablewidth}{0.5pt}}
\end{center}
\vspace{-4mm} \caption{\label{tab2} Accuracy test for the isotropic
vortex propagation problem.}
\end{table}

\subsection{One dimensional Riemann problem}
The first one is Sod problem \cite{Case-Sod}, the computational
domain is $[0,1]$ and the ratio of specific heats takes
$\gamma=1.4$. The initial condition is given by
\begin{equation*}
(\rho,u,p)=\left\{\begin{aligned}
&(1, 0, 1), 0<x<0.5,\\
&(0.125,0,0.1),  0.5<x<1.
\end{aligned} \right.
\end{equation*}
The density, velocity and pressure distributions with 100 meshes and
the exact solution at $t=0.2$ are given in Fig.\ref{1d}, and the
numerical results agree well with the exact solutions.

The second one is the Woodward-Colella blast wave problem
\cite{Case-Woodward}. The computational domain is $[0,100]$ with 400 mesh points and with reflected boundary condition on both ends.
The ratio of specific
heats also takes $\gamma=1.4$. The initial condition are given as
follows
\begin{equation*}
(\rho,u,p)=\left\{\begin{aligned}
&(1, 0, 1000), 0\leq x<10,\\
&(1, 0, 0.01), 10\leq x<90,\\
&(1, 0, 100),  90\leq x\leq 100.
\end{aligned} \right.
\end{equation*}
The density, velocity and pressure distributions at $t=3.8$ are
presented in Fig \ref{1d}, which are compared with the reference
solutions obtained by the second-order BGK scheme with van Leer
limiter. The figures show that the scheme can well resolve the
strong shock and contact discontinuities, particularly for the local
extreme values.

\begin{figure}[!h]
\centering
\includegraphics[width=0.42\textwidth]{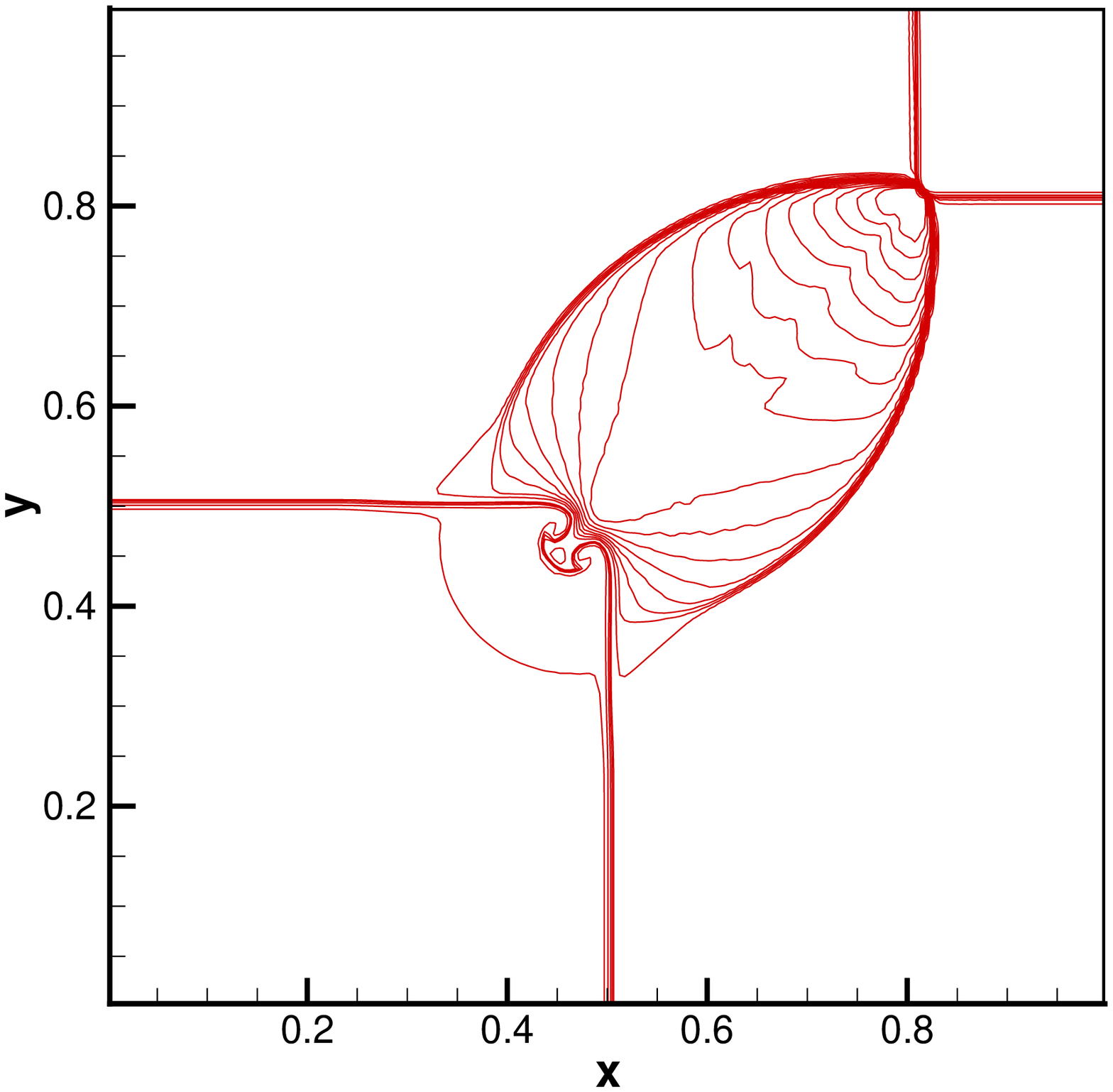}\includegraphics[width=0.42\textwidth]{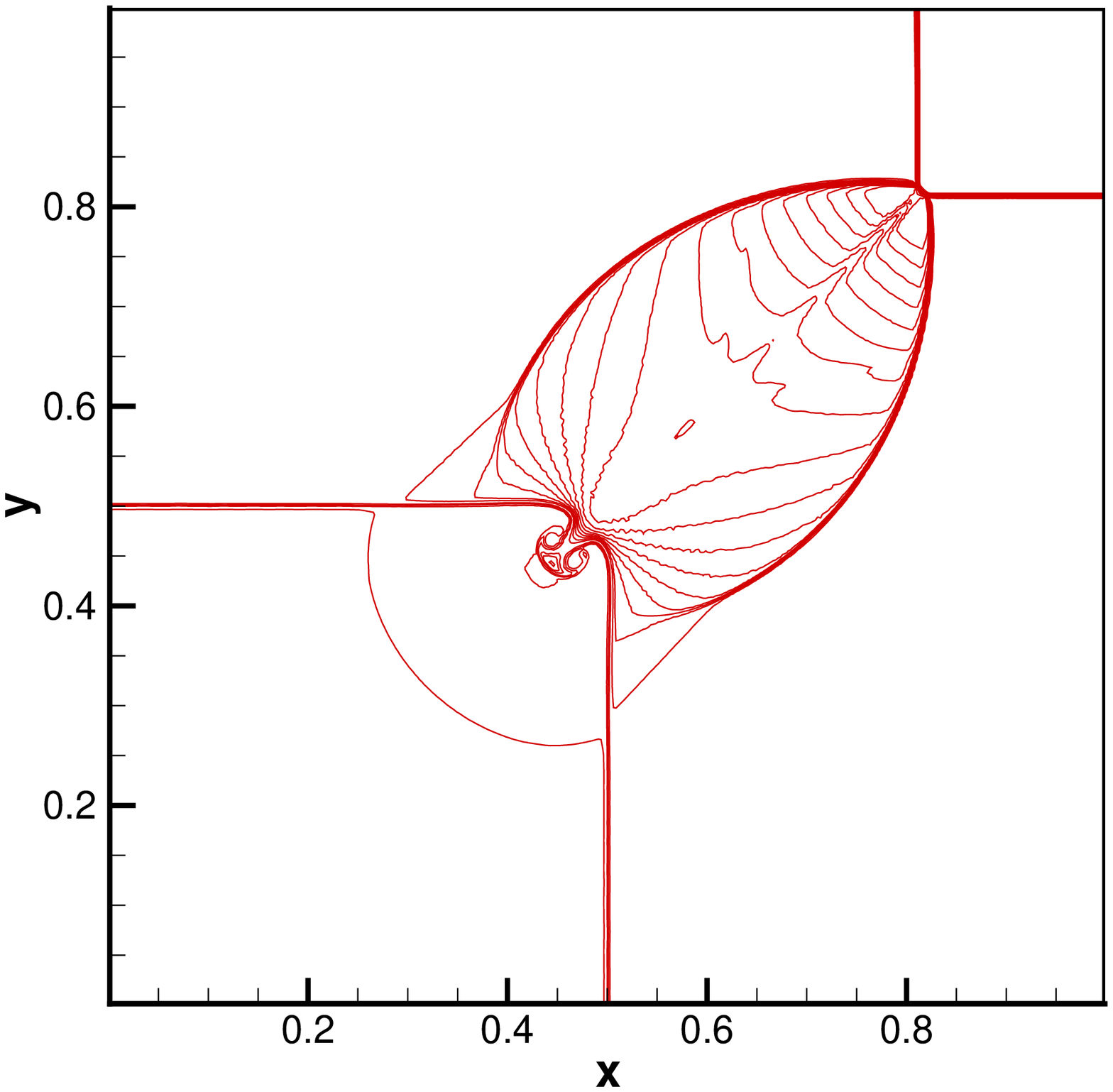}
\includegraphics[width=0.42\textwidth]{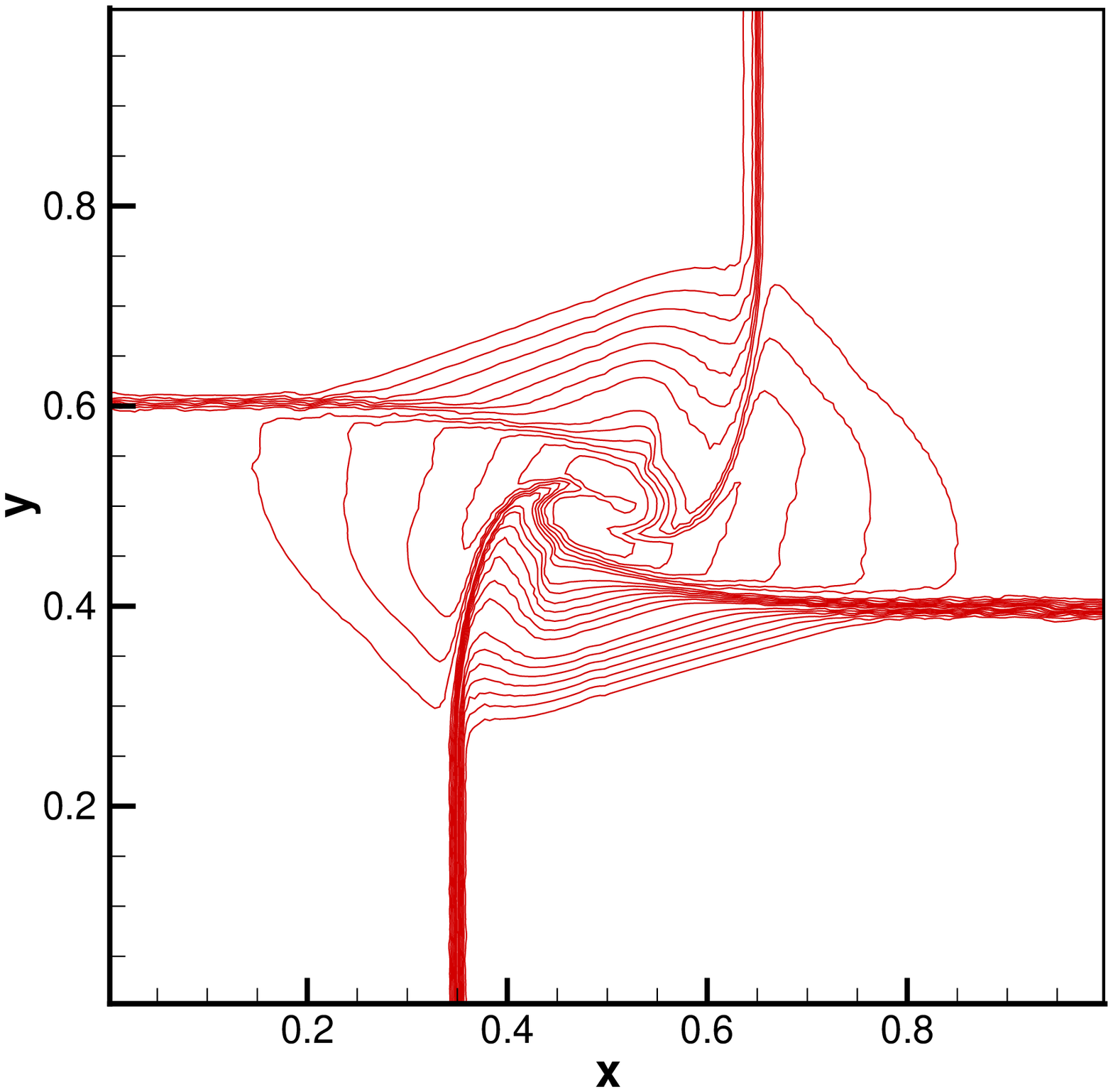}\includegraphics[width=0.42\textwidth]{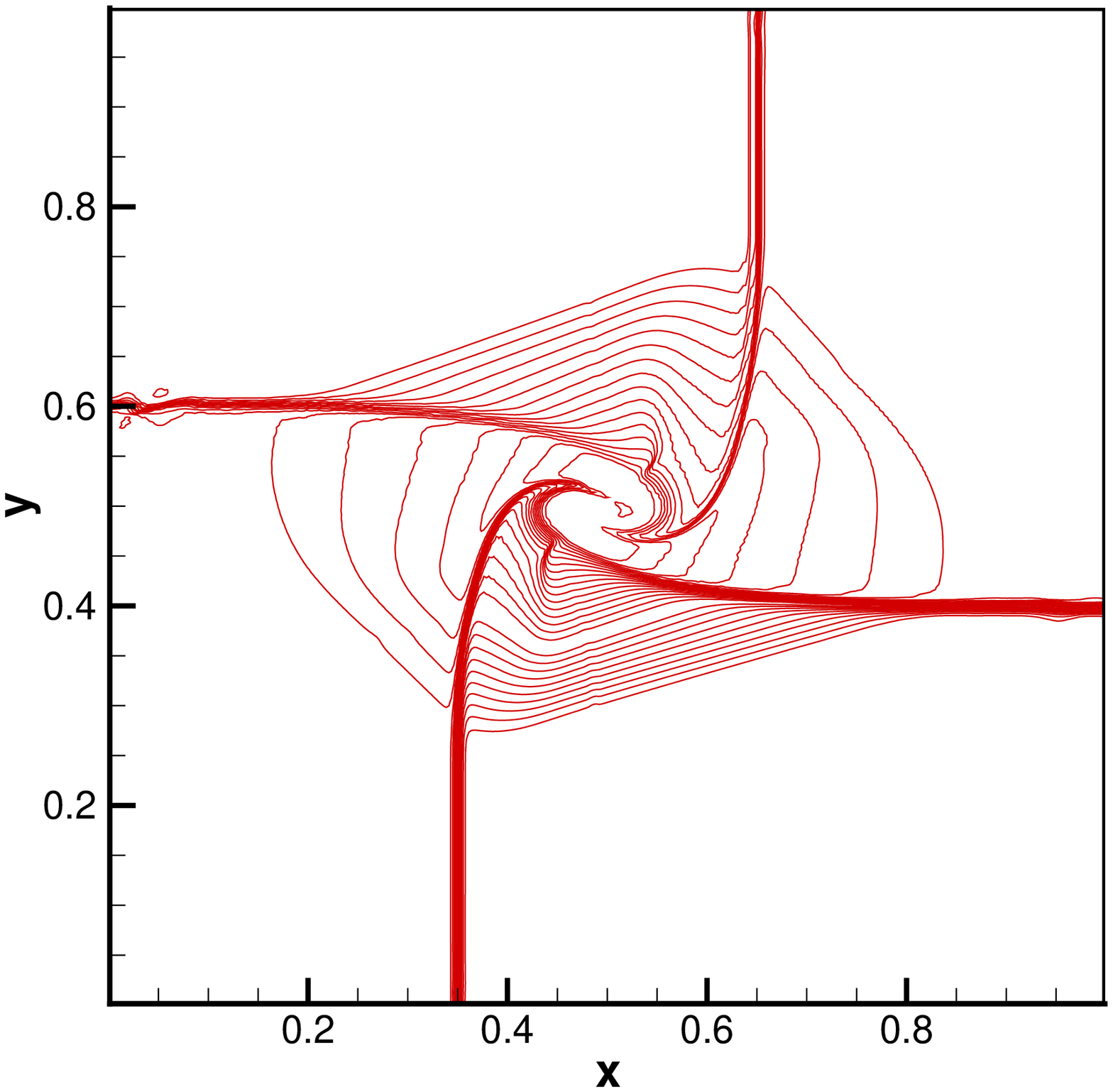}
\caption{\label{2d-riemann} Two-dimensional Riemann problem for case
1 (top) and case 2 (bottom): the density distributions at $t=0.2$.
Left: 200$\times$200 cells; right: 400$\times$400 cells. }
\end{figure}

\subsection{Two dimensional Riemann problem}
In this case, two 2-D Riemann problems are tested to verify the
capability in capturing complex two dimensional wave configurations.
Both cases were presented in \cite{Case-Liu}. The
computational domain is $[0,1]\times[0,1]$ and $\gamma=1.4$. The
initial condition for the first case is
\begin{equation*}
(\rho,u,v,p)=\left\{\begin{aligned}
         &(0.5313,0,0,0.4), &x>0.5,y>0.5,\\
         &(1,0.7276,0,1), &x<0.5,y>0.5,\\
         &(0.8,0,0,1), &x<0.5,y<0.5,\\
         &(1,0,0.7276,1), &x>0.5,y<0.5.
                          \end{aligned} \right.
                          \end{equation*}
The initial condition for the second case is
\begin{equation*}
(\rho,u,v,p)=\left\{\begin{aligned}
         &(1 ,0.75,-0.5,1), &x>0.5,y>0.5,\\
         &(2,0.75,0.5,1), &x<0.5,y>0.5,\\
         &(1,-0.75,0.5,1), &x<0.5,y<0.5,\\
         &(3.,-0.75,-0.5,1), &x>0.5,y<0.5.
                          \end{aligned} \right.
                          \end{equation*}
Non-reflecting boundary conditions are used in $x$ and $y$
directions in the computation. The density distributions at $t=0.2$
shown in Fig.\ref{2d-riemann} for these cases with $200\times200$
and $400\times400$ cells, and the results show that the current
scheme well resolve the flow structure.

\begin{figure}[!h]
\centering
\includegraphics[width=0.58\textwidth]{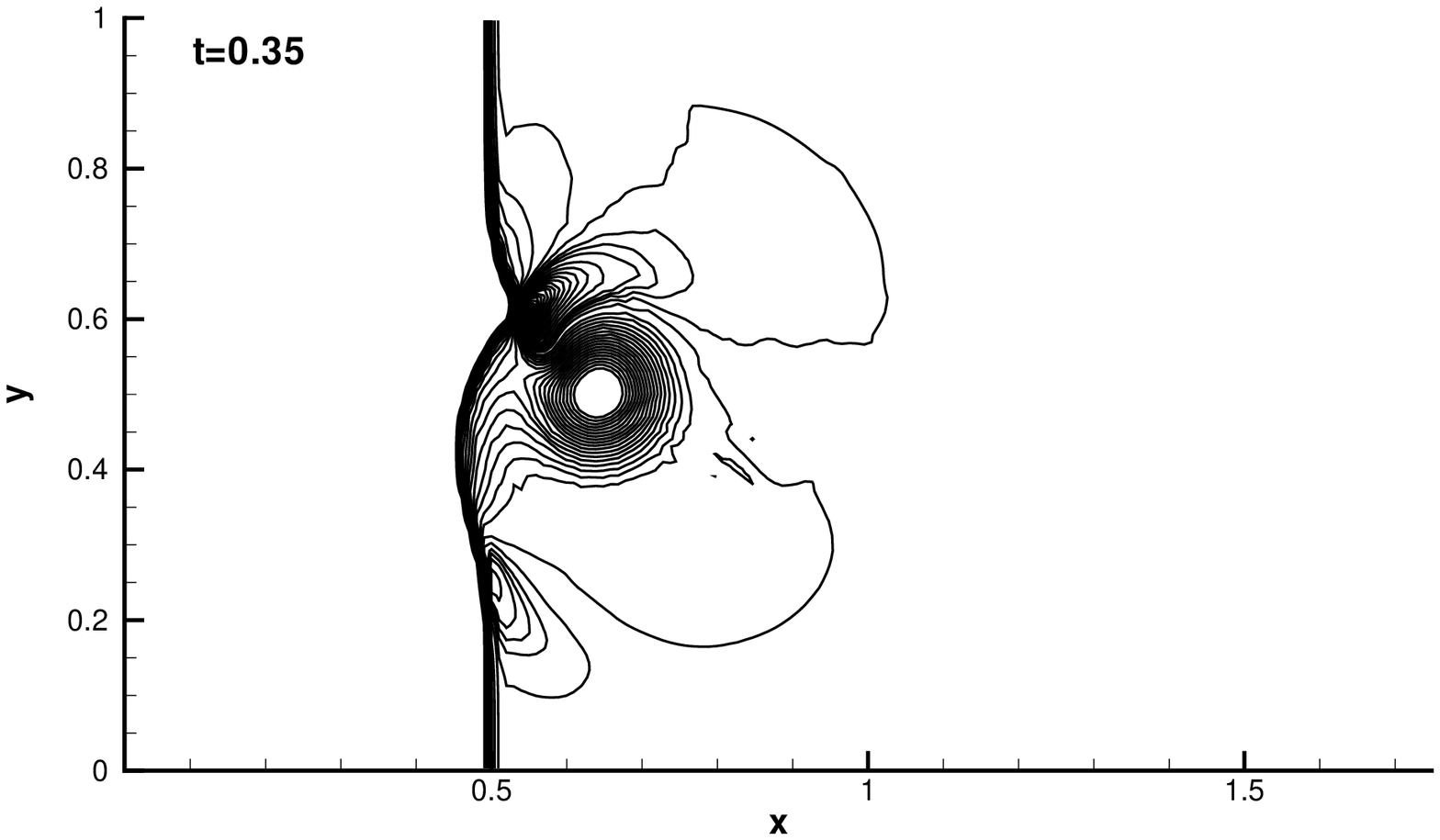}
\includegraphics[width=0.58\textwidth]{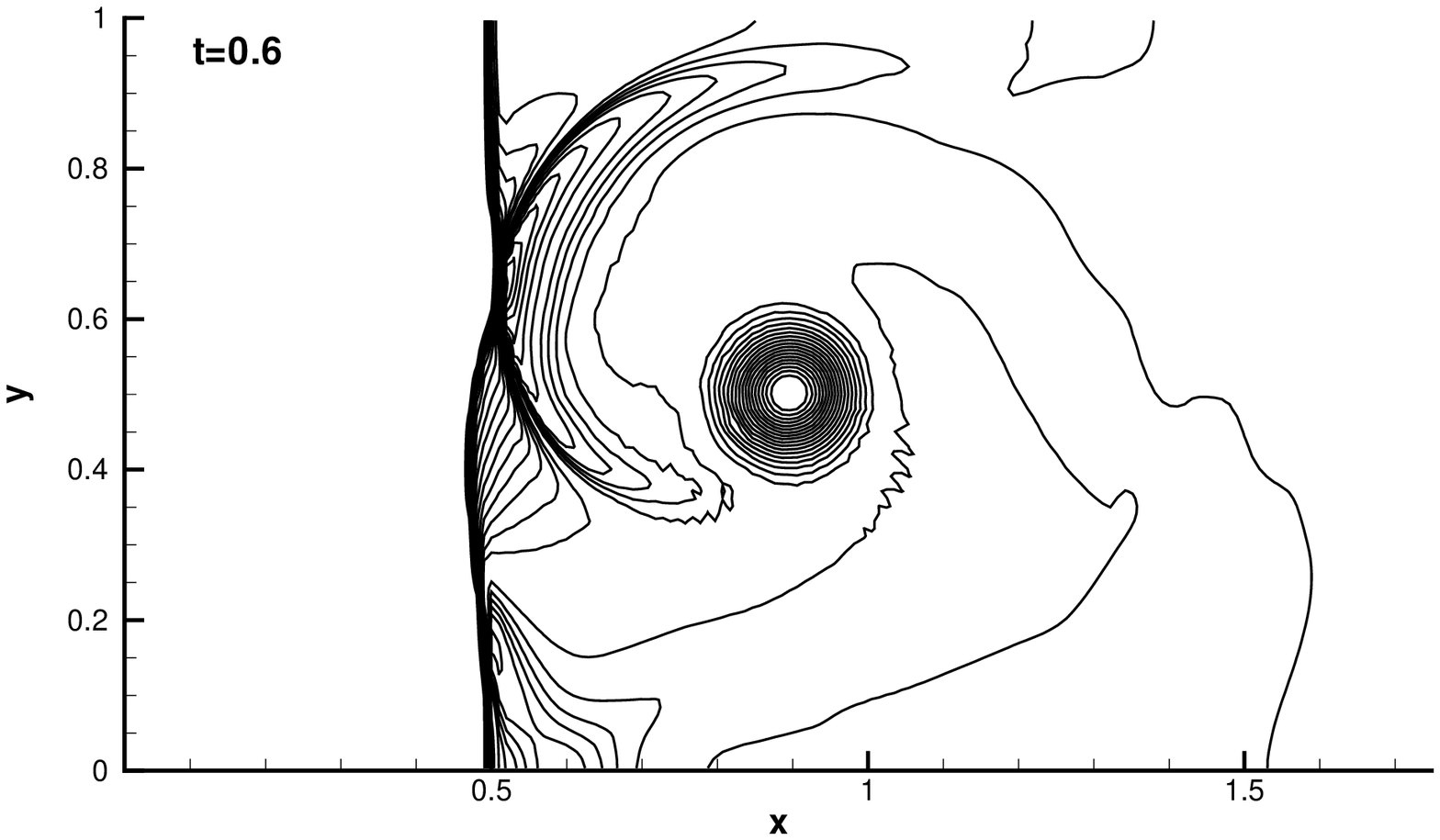}
\includegraphics[width=0.58\textwidth]{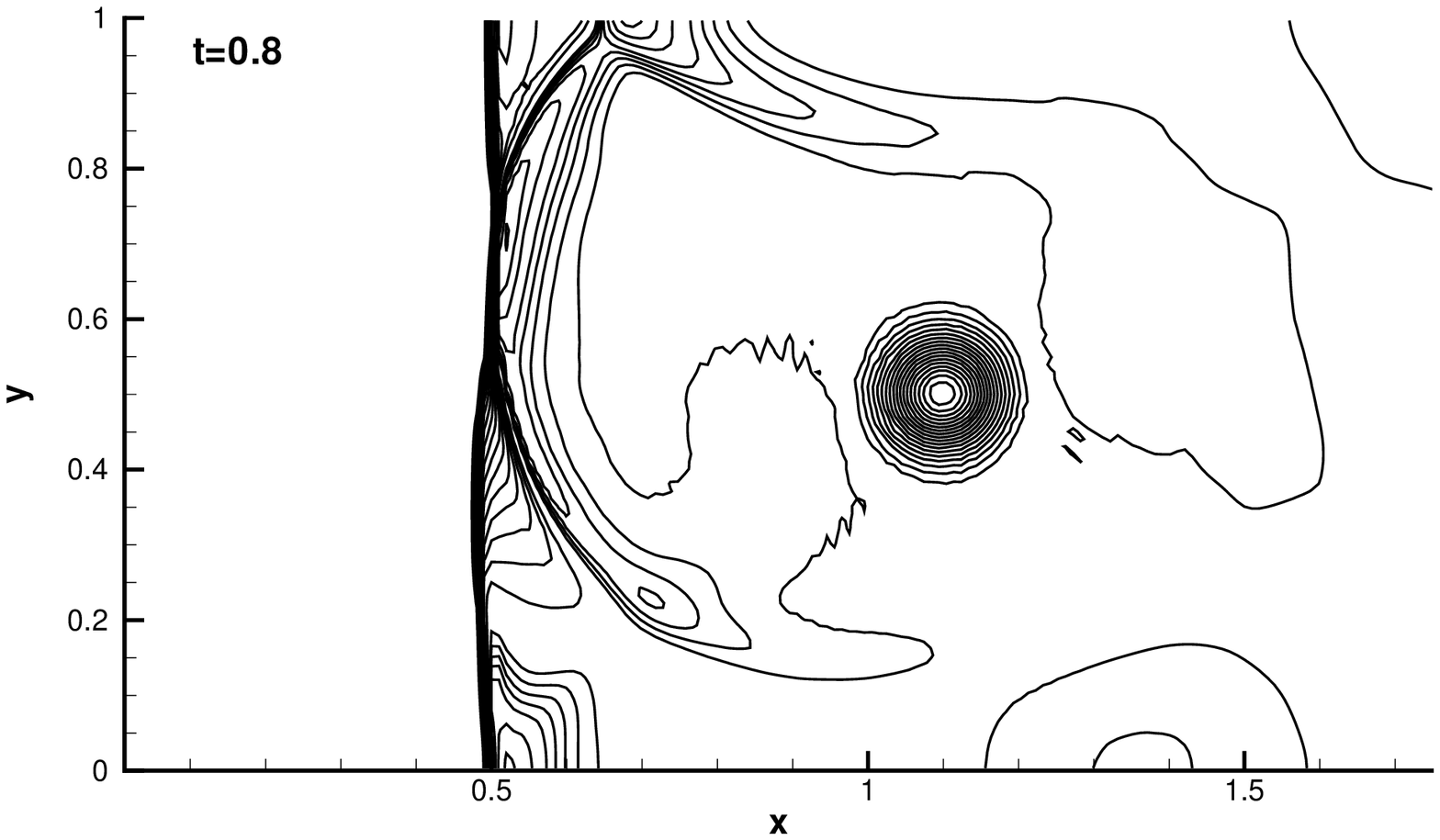}
\caption{\label{shock-vortex1} The pressure distribution for
two-dimensional shock vortex interaction at $t=0.35, 0.6$ and $0.8$
with $201\times101$ uniform mesh points. }
\end{figure}

\subsection{Shock vortex interaction}
\begin{figure}[!h]
\centering
\includegraphics[width=0.5\textwidth]{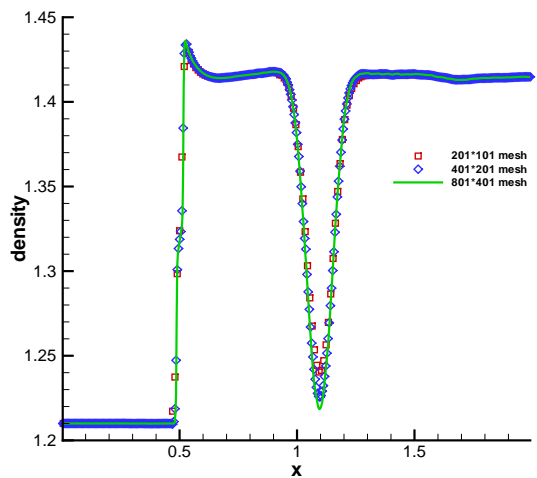}
\caption{\label{shock-vortex2} The density distribution for
two-dimensional shock vortex interaction at $t=0.8$ along the
horizontal symmetric line $y= 0.5$ with $201\times101$,
$401\times201$ and $801\times401$ uniform mesh points. }
\end{figure}

This model problem describes the interaction between a stationary
shock and a vortex for the inviscid flow, which was presented in
\cite{WENO2}. The computational domain is taken to be $[0,
2]\times[0, 1]$. A stationary Mach $1.1$ shock is positioned at
$x=0.5$ and normal to the $x$-axis. The left upstream state is
$(\rho, u, v, p) = (Ma^2,\sqrt{\gamma}, 0, 1)$, where $\gamma=1.4$
is the specific heat ratio and $Ma$ is the Mach number. A small
vortex is obtained through a perturbation on the mean flow with the
velocity $(u, v)$, temperature $T=p/\rho$ and entropy
$S=\ln(p/\rho^\gamma)$, and the perturbation is expressed as
\begin{align*}
&(\delta u,\delta v)=\kappa\eta e^{\mu(1-\eta^2)}(\sin\theta,-\cos\theta),\\
&\delta
T=-\frac{(\gamma-1)\kappa^2}{4\mu\gamma}e^{2\mu(1-\eta^2)},\delta
S=0,
\end{align*}
where $\eta=r/r_c$, $r=\sqrt{(x-x_c)^2+(y-y_c)^2}$, $(x_c,
y_c)=(0.25, 0.5)$ is the center of the vortex. $\kappa$ indicates
the strength of the vortex, $\mu$ controls the decay rate of the
vortex  and $r_c$ is the critical radius for which the vortex has
the maximum strength. In the computation, $\kappa=0.3$, $\mu=0.204$
and $r_c=0.05$. The reflected boundary condition is used on the top
and bottom boundaries. The pressure distributions with
$201\times101$ mesh points at $t=0.35$, $0.6$ and $0.8$ are shown in
Fig.\ref{shock-vortex1}. By $t=0.8$, one branch of the shock
bifurcations has reached the top boundary and been reflected, and
the reflection is well captured. The detailed density distribution
along the center horizontal line with $201\times101$, $401\times201$
and $801\times401$ mesh points at $t=0.8$ are shown in
Fig.\ref{shock-vortex2}.

\begin{figure}[!h]
\centering
\includegraphics[width=0.8\textwidth]{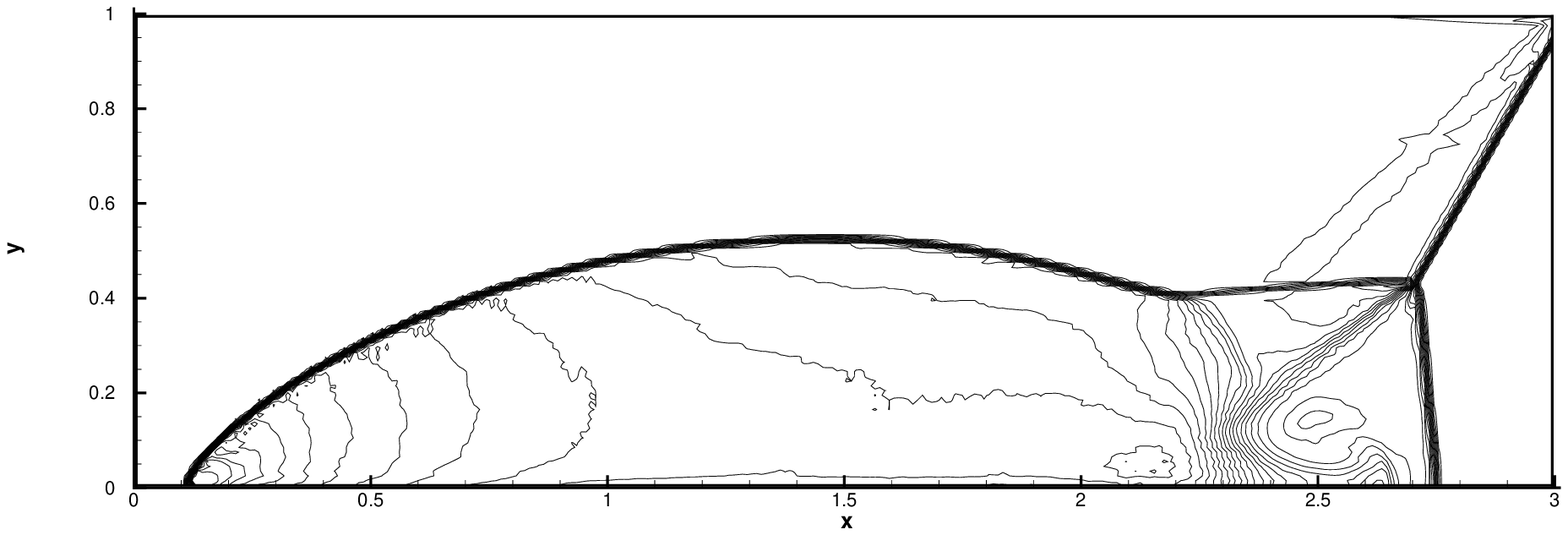}
\includegraphics[width=0.8\textwidth]{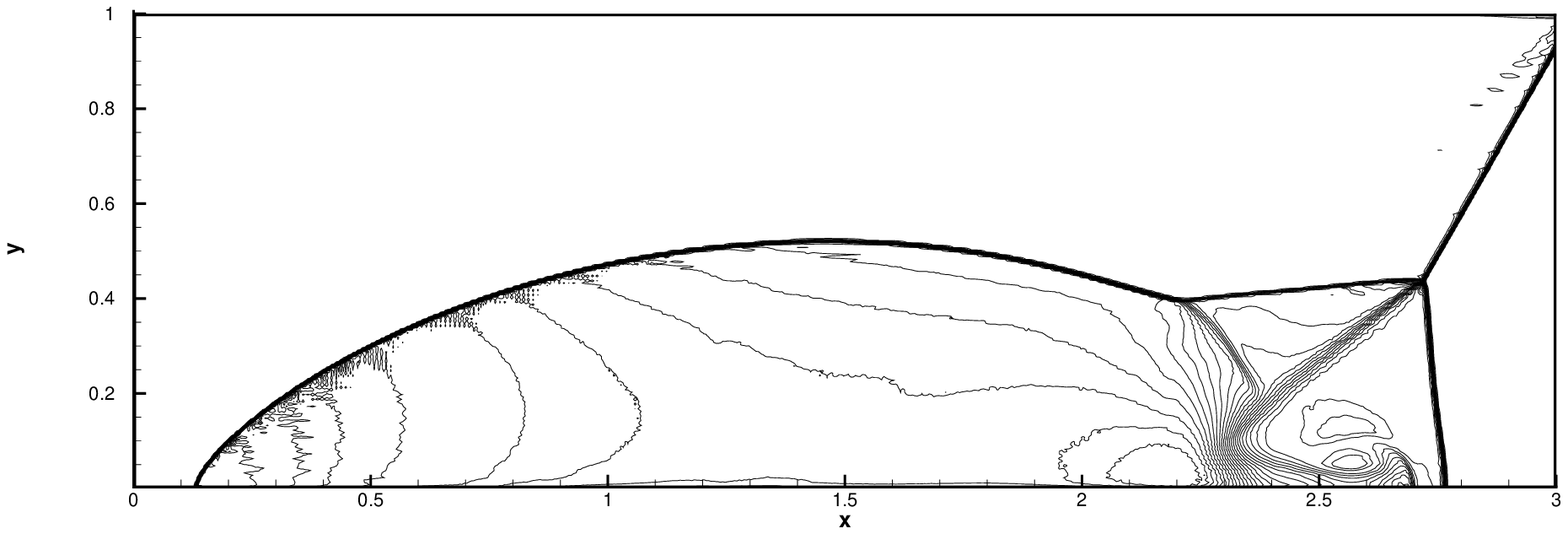}
\caption{\label{double-mach-1} The density distributions of double
mach reflection problem with the cell size $\Delta x=\Delta y=1/100$
(top) and $\Delta x=\Delta y=1/200$ (bottom) at $t=0.2$.}
\end{figure}

\subsection{Double Mach reflection problem}
This problem was extensively studied by Woodward and Colella
\cite{Case-Woodward} for the inviscid flow. The computational domain
is $[0,4]\times[0,1]$, and a solid wall lies at the bottom of the
computational domain starting from $x =1/6$. Initially a
right-moving Mach 10 shock is positioned at $(x,y)=(1/6, 0)$, and
made a $60^\circ$ angle with the x-axis. The initial pre-shock and
post-shock conditions are
\begin{align*}
(\rho, U, V, p)&=(8, 4.125\sqrt{3}, -4.125,
116.5),\\
(\rho, U, V, p)&=(1.4, 0, 0, 1).
\end{align*}
The reflective boundary condition is used at the wall, while for the
rest of bottom boundary, the exact post-shock condition is imposed.
At the top boundary, the flow values are set to describe the exact
motion of the Mach 10 shock. The density distributions with
$400\times100$ and $800\times200$ mesh points at $t=0.2$ are shown
in Fig.\ref{double-mach-1}. The current compact scheme resolves the
flow structure under the triple Mach stem clearly.

\begin{figure}[!h]
\centering
\includegraphics[width=0.725\textwidth]{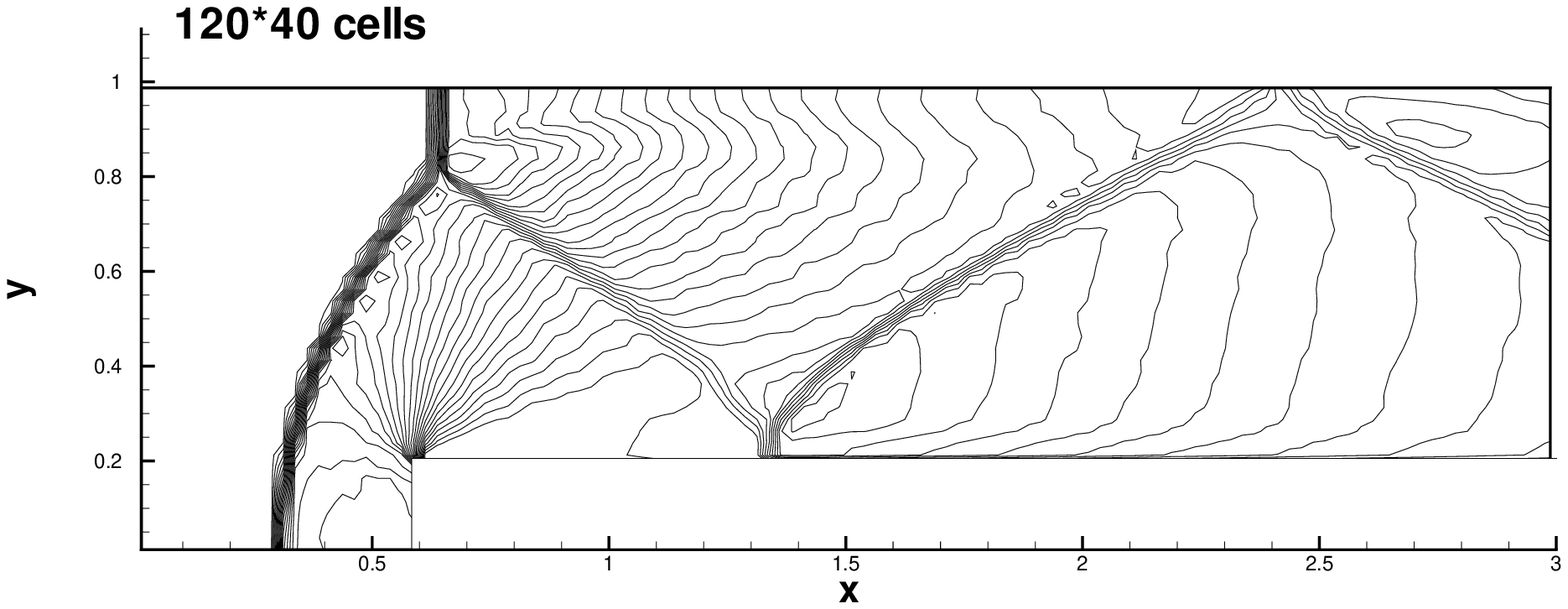}
\includegraphics[width=0.725\textwidth]{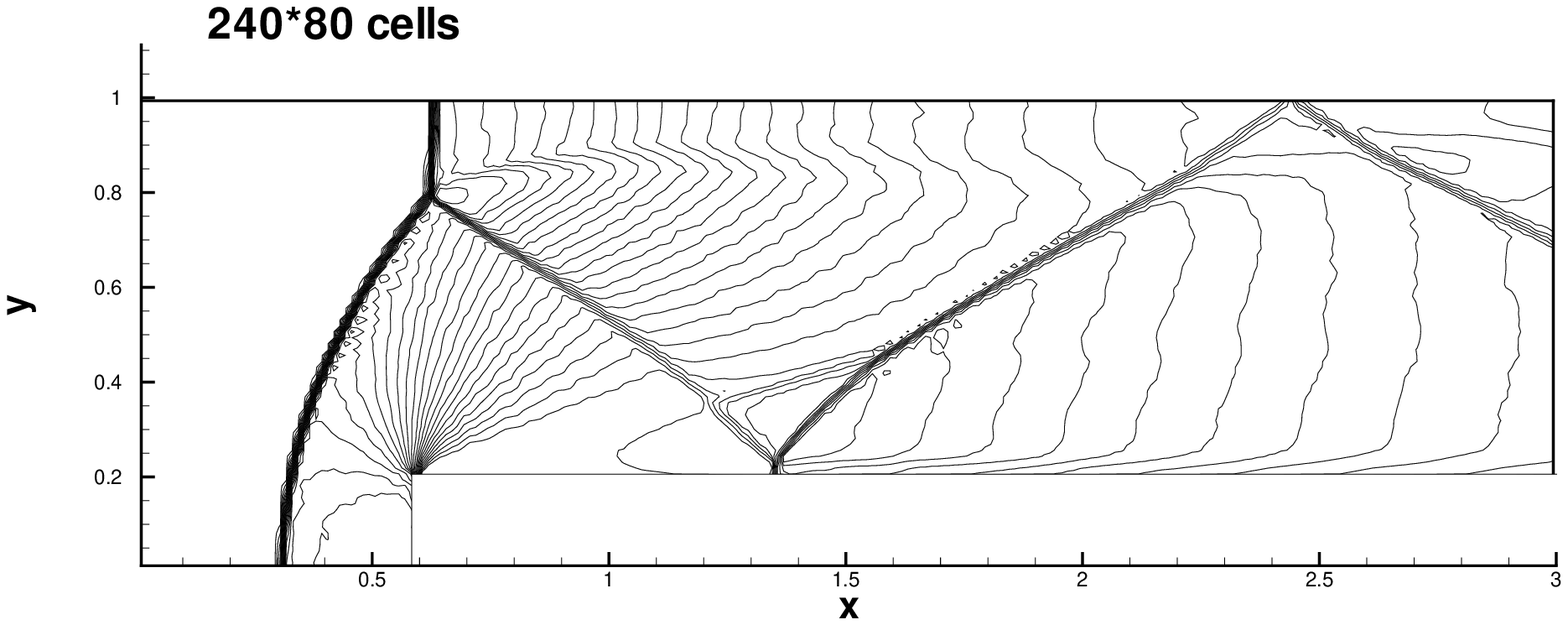}
\includegraphics[width=0.725\textwidth]{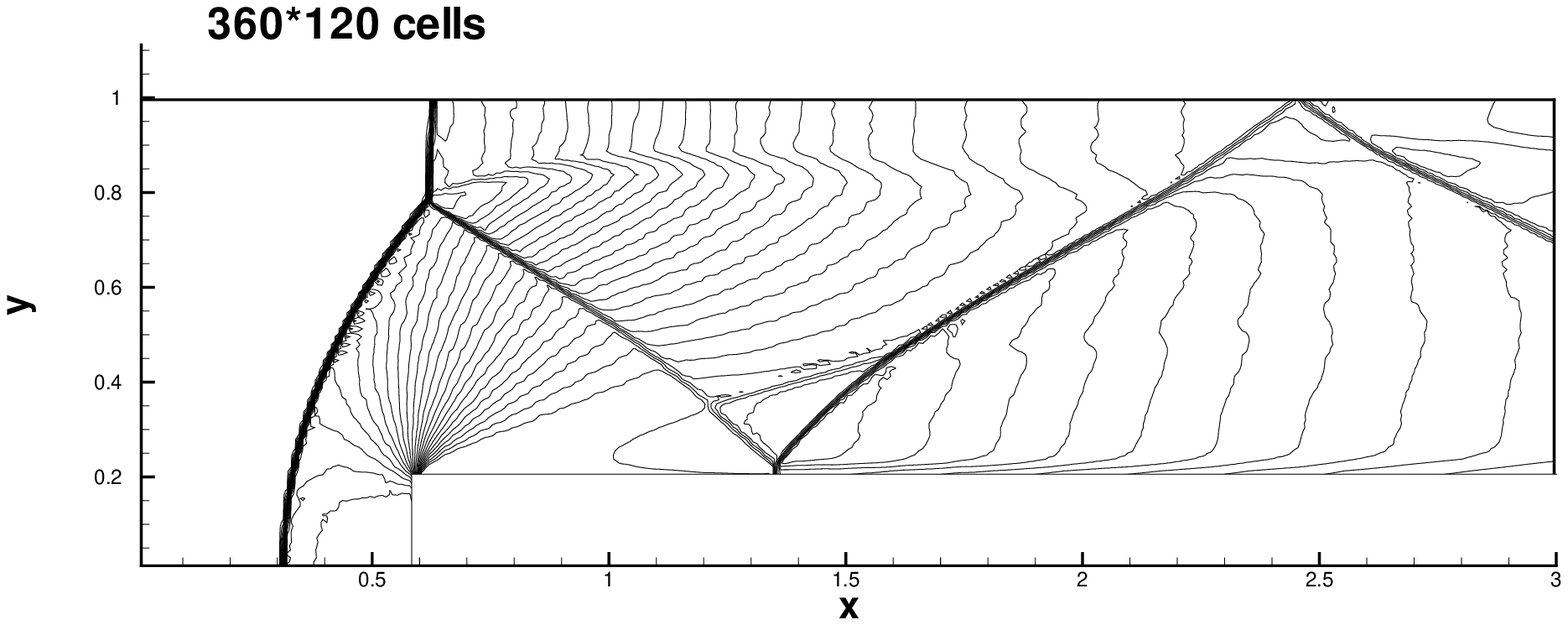}
\caption{\label{front-step-1} The density distribution of the front
step problem with $120\times40$, $240\times80$ and $360\times120$
mesh points for the inviscid flow at $t=4$.}
\end{figure}

\subsection{Front step problem}
The front step problem was again studied extensively by Woodward and
Colella \cite{Case-Woodward} for the inviscid flow. The
computational domain $[0,3]\times[0,1]$. The step is located at
$x=0.6$ with height $0.2$ in the tunnel. Initially, a right-going
Mach 3 flow is used. Reflective boundary conditions are used along
the walls of the tunnel, and inflow and outflow boundary conditions
are used at the entrance and the exit. The corner of the step is
center of a rarefaction fan, hence it is a singularity point. With
nothing special done at this point, the flow will be affected by the
erroneous entropy layer. To minimize the numerical error generated
at the corner of the step, the flow variables around the corner are
modified according to \cite{Case-Woodward} in the computation. The
density distributions with $120\times40$, $240\times80$ and
$360\times120$ mesh points are presented in Fig.\ref{front-step-1}
at $t=4$.

\subsection{Viscous shock tube problem}
This problem was introduced in \cite{Case-Daru,Case-Sjogreen} to
test the performances of different schemes for viscous flows. In
this case, an ideal gas is at rest in a two-dimensional unit box
$[0,1]\times[0,1]$. A membrane is located at $x=0.5$ separates two
different states of the gas and the dimensionless initial states are
\begin{equation*}
(\rho,u,p)=\left\{\begin{aligned}
&(120, 0, 120/\gamma), 0<x<0.5,\\
&(1.2, 0, 1.2/\gamma),  0.5<x<1,
\end{aligned} \right.
\end{equation*}
where $\gamma=1.4$ and Prandtl number $Pr=0.73$.

\begin{figure}[!h]
\centering
\includegraphics[width=0.88\textwidth]{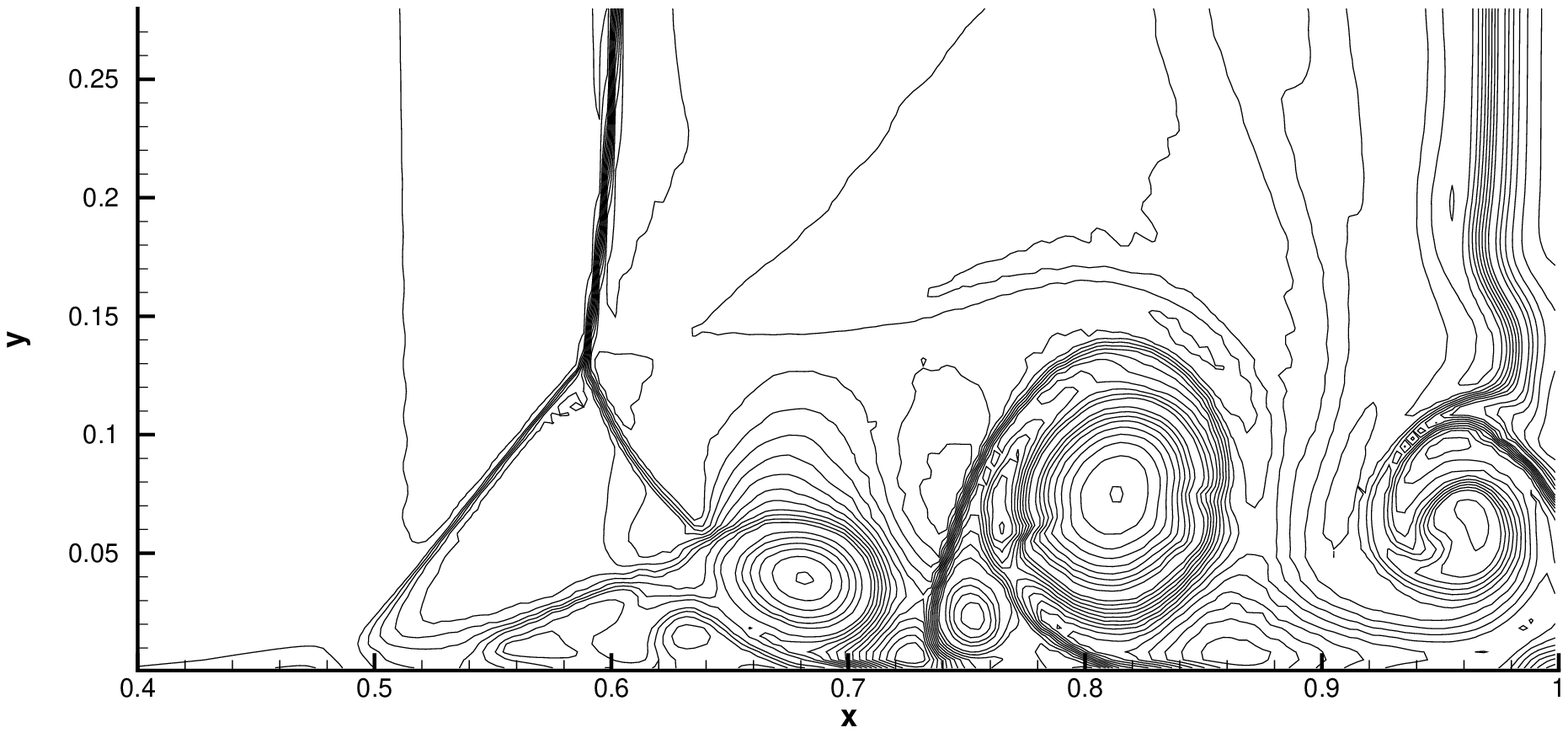}
\caption{\label{shock-boundary1} Reflected shock-boundary layer
interaction. The density distribution at $t=1$ with $300\times150$
mesh points with $Re=200$.}
\end{figure}

\begin{figure}[!h]
\centering
\includegraphics[width=0.88\textwidth]{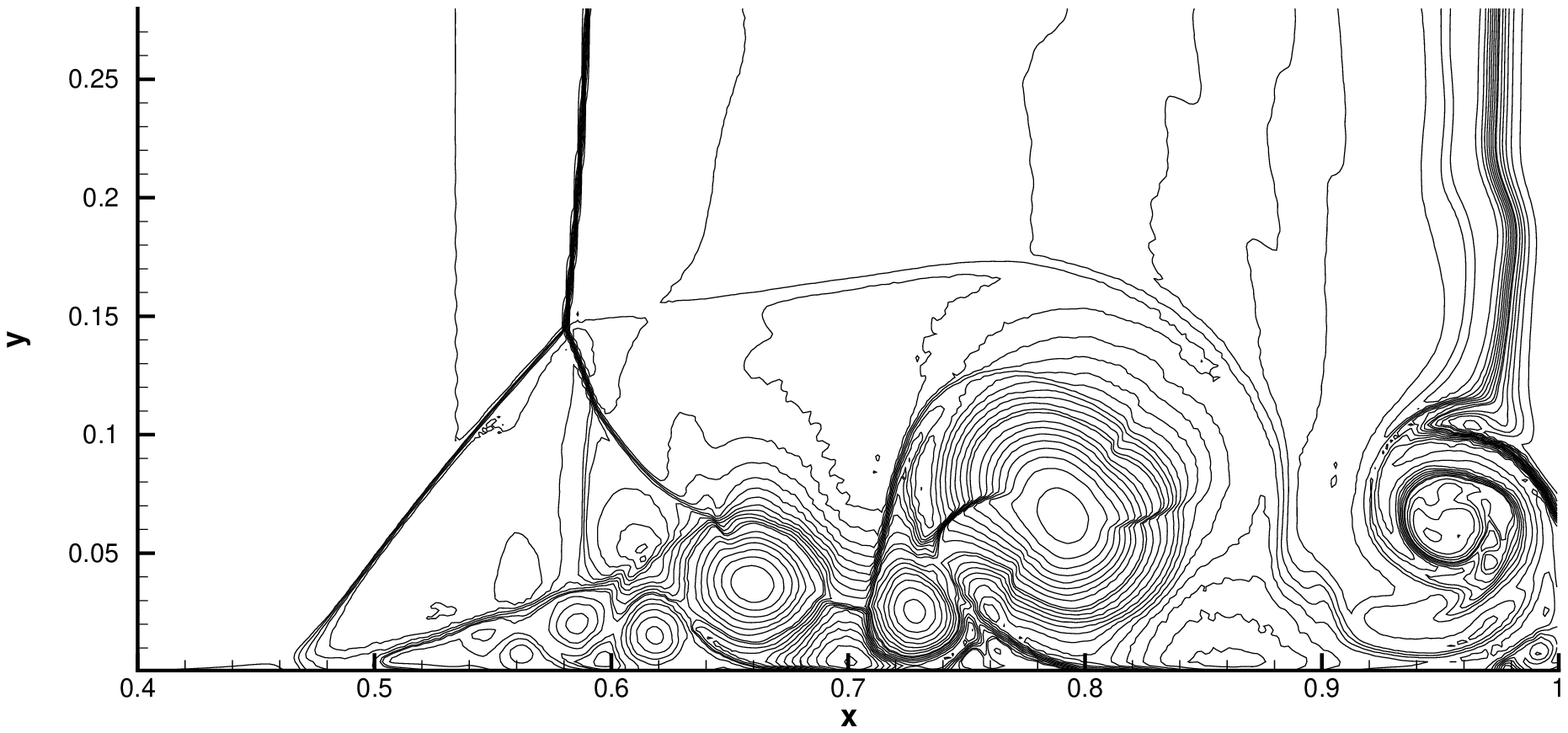}
\includegraphics[width=0.88\textwidth]{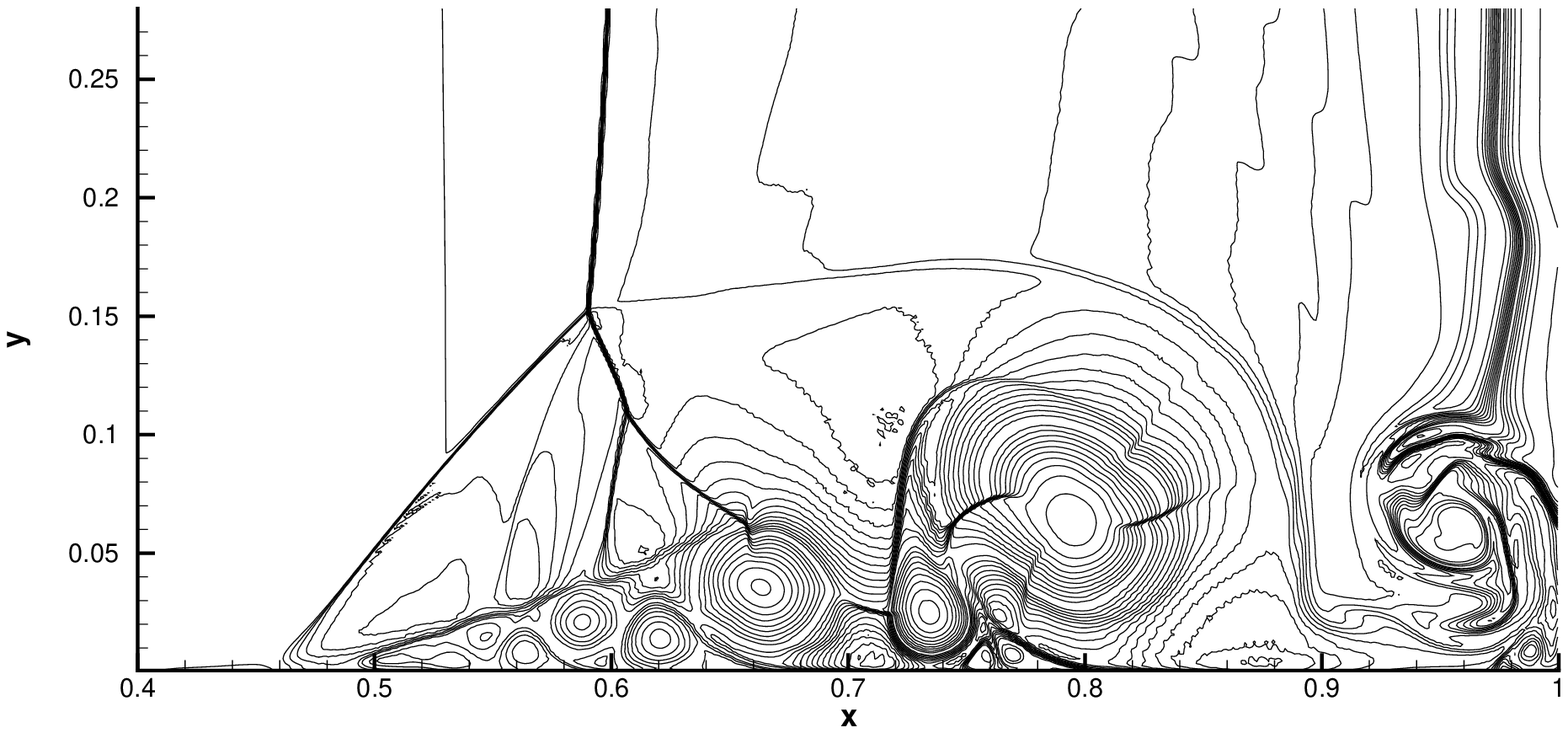}
\caption{\label{shock-boundary2} Reflected shock-boundary layer
interaction. The density distribution at $t=1$ with $600\times300$
(top) and $1000\times500$ (bottom) mesh points with $Re=1000$.}
\centering
\includegraphics[width=0.65\textwidth]{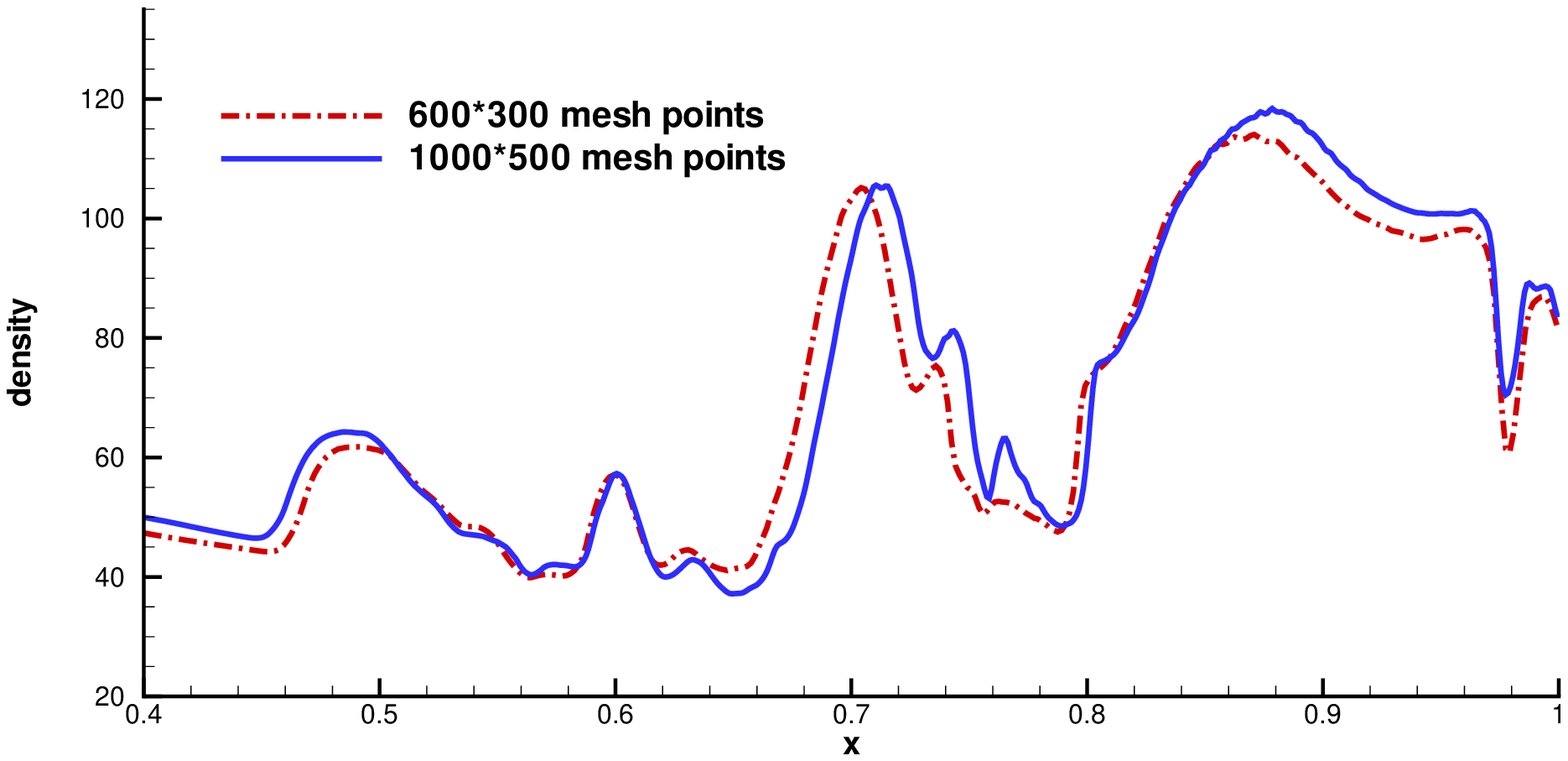}
\caption{\label{shock-boundary-x} Reflected shock-boundary layer
interaction. The density distribution at $t=1$ with along the lower
wall with $600\times300, 1000\times500$ mesh points with $Re=1000$.}
\end{figure}

The membrane is removed at time zero and wave interaction occurs. A
shock wave, followed by a contact discontinuity, moves to the right
with Mach number $Ma=2.37$ and reflects at the right end wall. After
the reflection, it interacts with the contact discontinuity. The
contact discontinuity and shock wave interact with the horizontal
wall and create a thin boundary layer during their propagation. The
solution will develop complex two-dimensional
shock/shear/boundary-layer interactions. This case is tested in the
computational domain $[0, 1]\times[0, 0.5]$, a symmetrical condition
is used on the top boundary $x\in[0, 1], y=0.5$ and non-slip
boundary condition and adiabatic condition for temperature are
imposed at solid wall boundaries. The density distribution at $t=1$
with $300\times150$ mesh points is shown in
Fig.\ref{shock-boundary1}. The complexity of the flow structure
increases as the Reynolds number increases. The density distribution
at $t=1$ with $600\times300$ and $1000\times500$ mesh points with
$Re=1000$ are shown in Fig.\ref{shock-boundary2}. The current scheme
can well resolve the complex flow structure. The density profiles
along the lower wall on with $600\times300$ and $1000\times500$ mesh
points with $Re=1000$ are presented in Fig.\ref{shock-boundary-x}.

\subsection{Lid-driven cavity flow}
\begin{figure}[!h]
\centering
\includegraphics[width=0.45\textwidth]{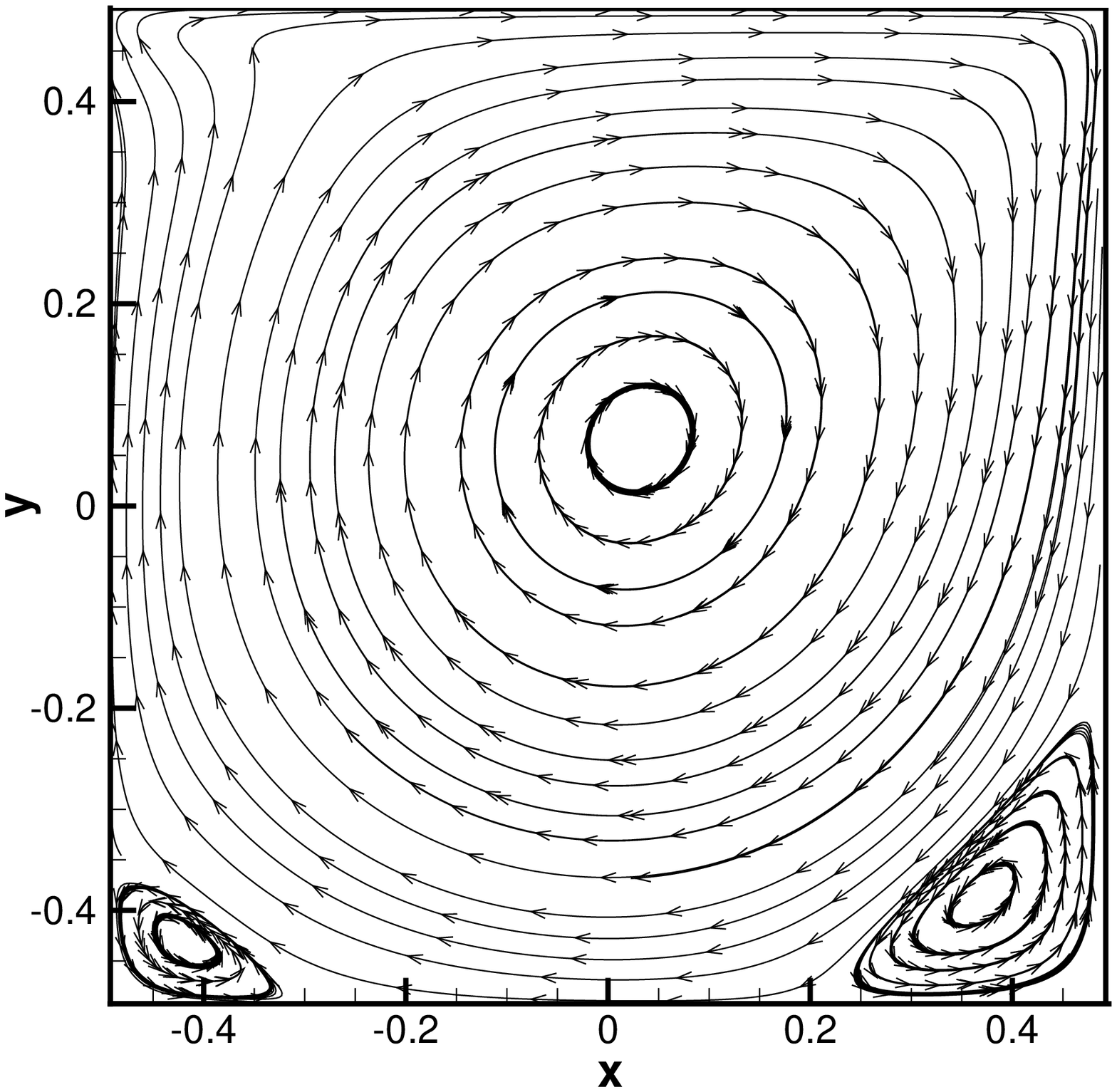}
\caption{\label{cavity-stream} The streamlines for the cavity flow
by $65\times65$ mesh points with $Re=1000$.}
\includegraphics[width=0.4\textwidth]{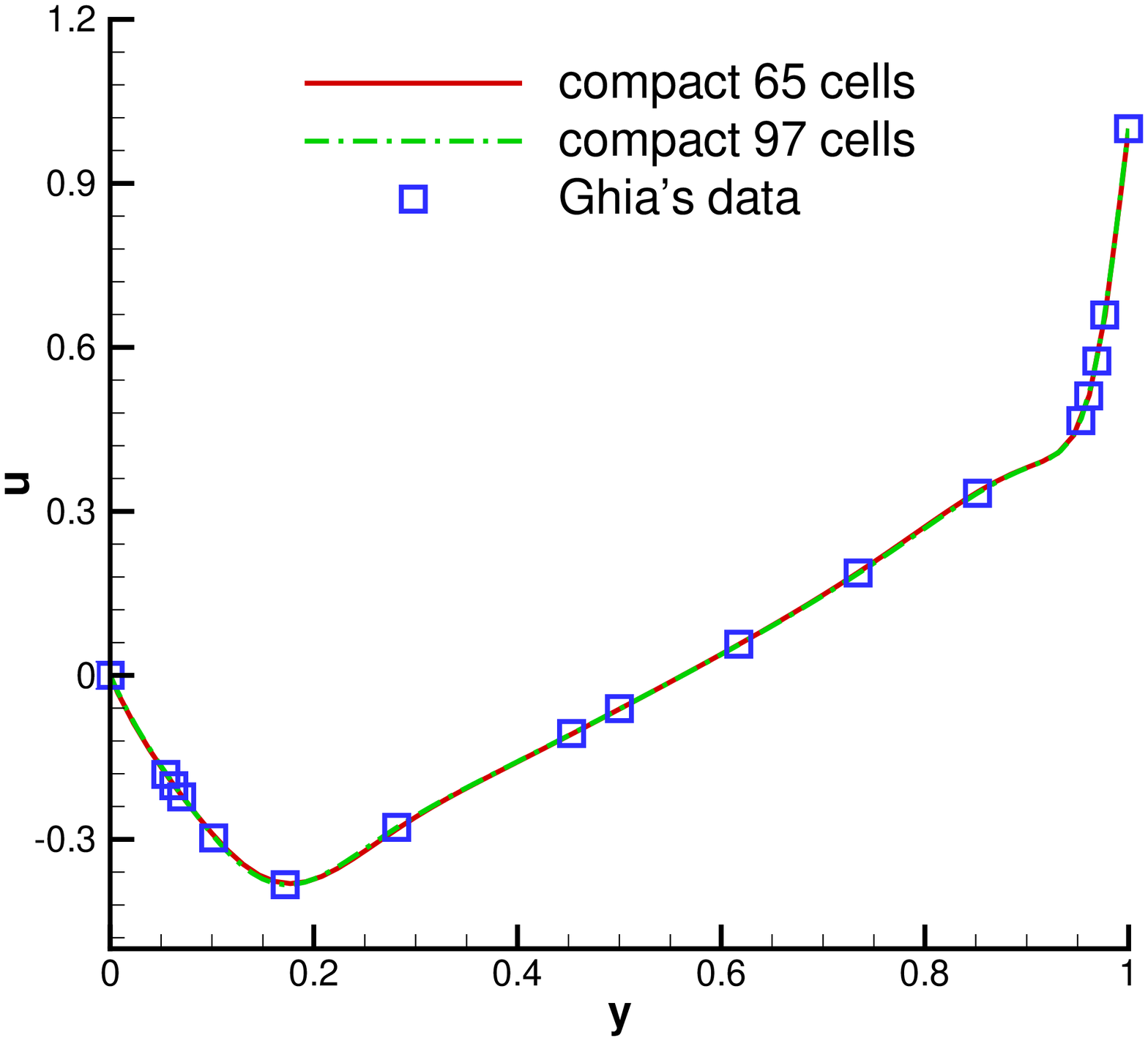}
\includegraphics[width=0.4\textwidth]{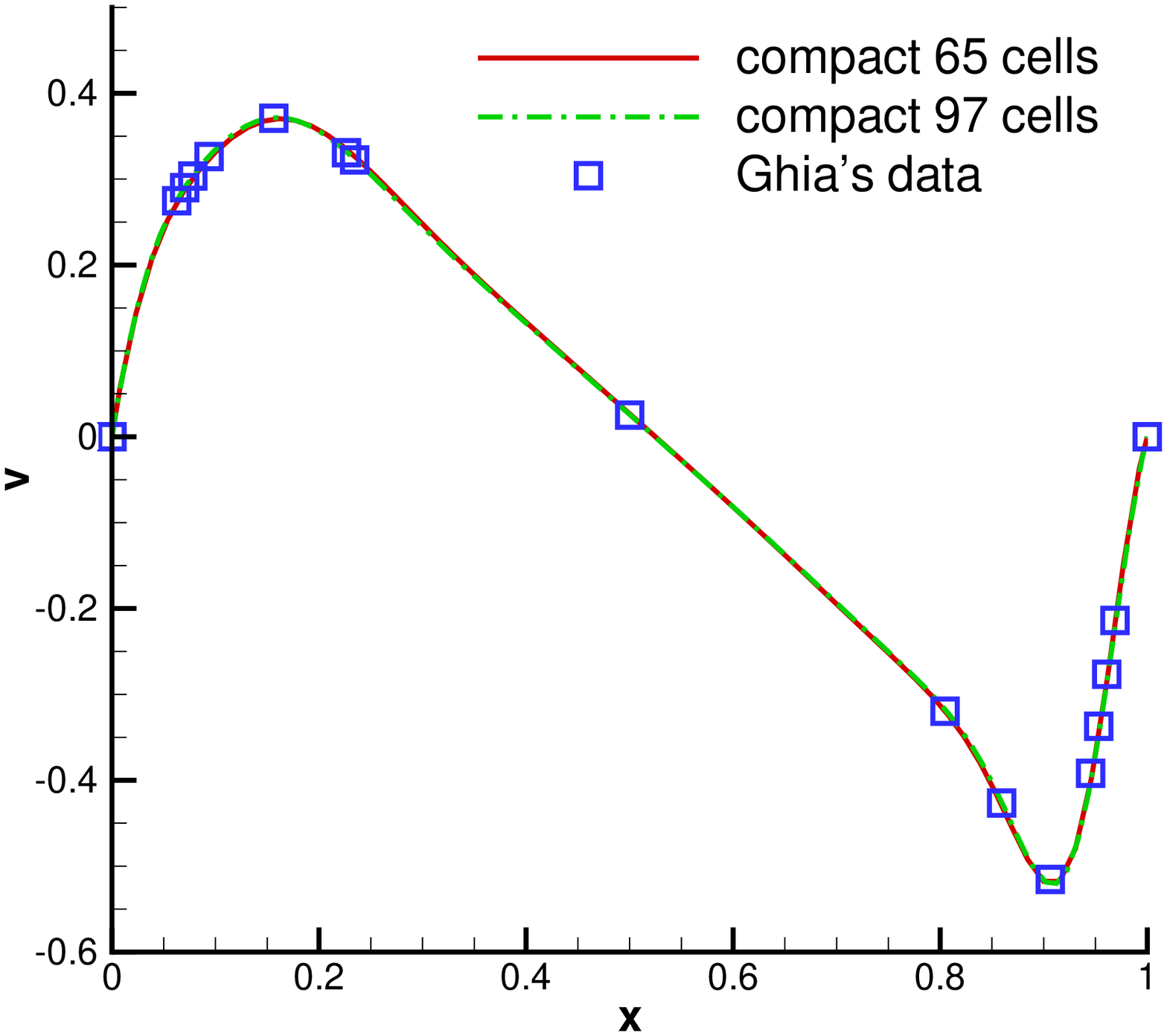}
\caption{\label{cavity-1} Lid-driven cavity flow: $U$-velocities
along vertical centerline line (left) and $V$-velocities along
horizontal centerline with $Re=1000$ with $65\times65$ and
$97\times97$ mesh points. The reference data is from
\cite{Case-Ghia}.}
\end{figure}

\begin{figure}[!h]
\centering
\includegraphics[width=0.4\textwidth]{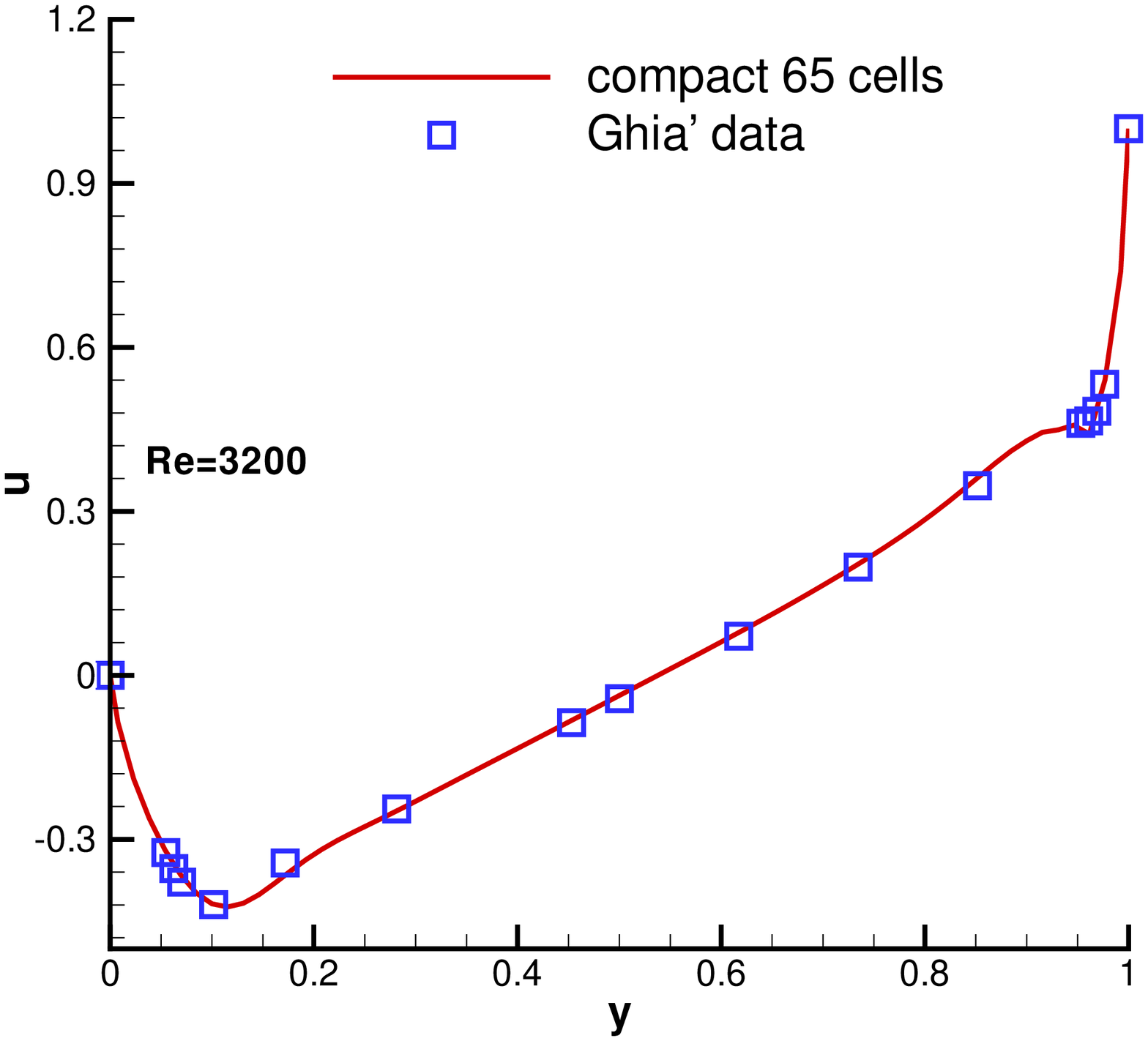}
\includegraphics[width=0.4\textwidth]{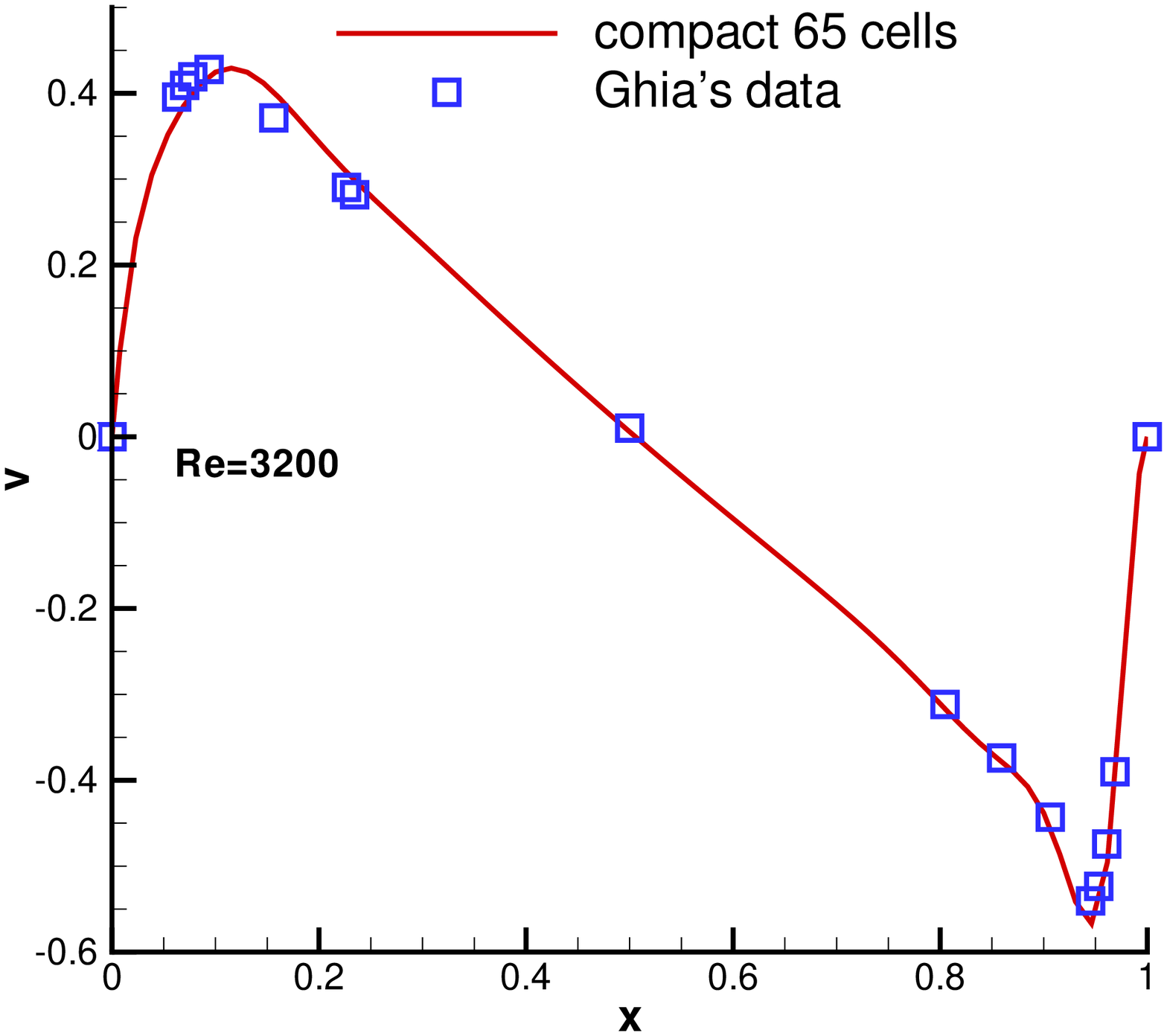}
\includegraphics[width=0.4\textwidth]{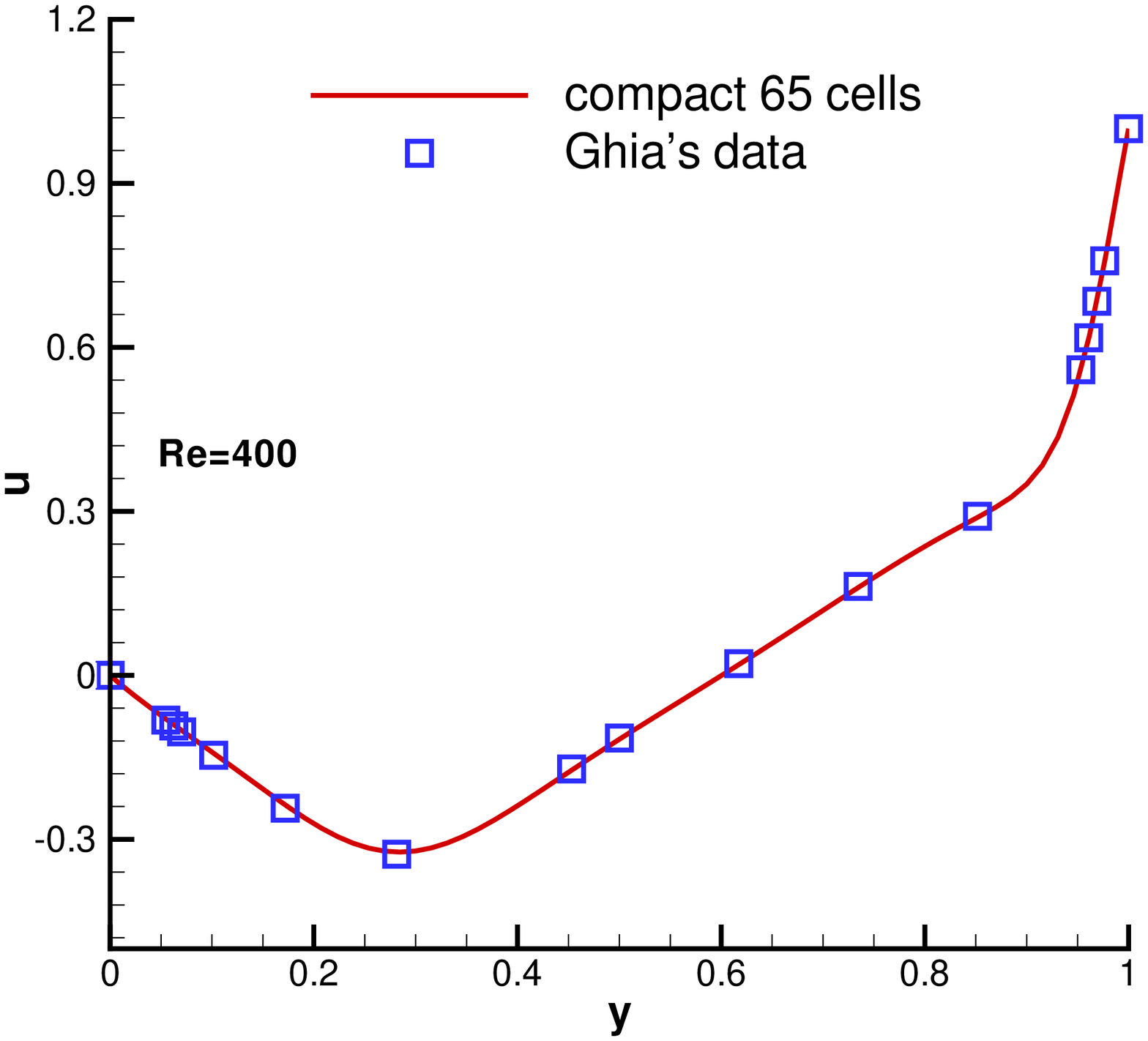}
\includegraphics[width=0.4\textwidth]{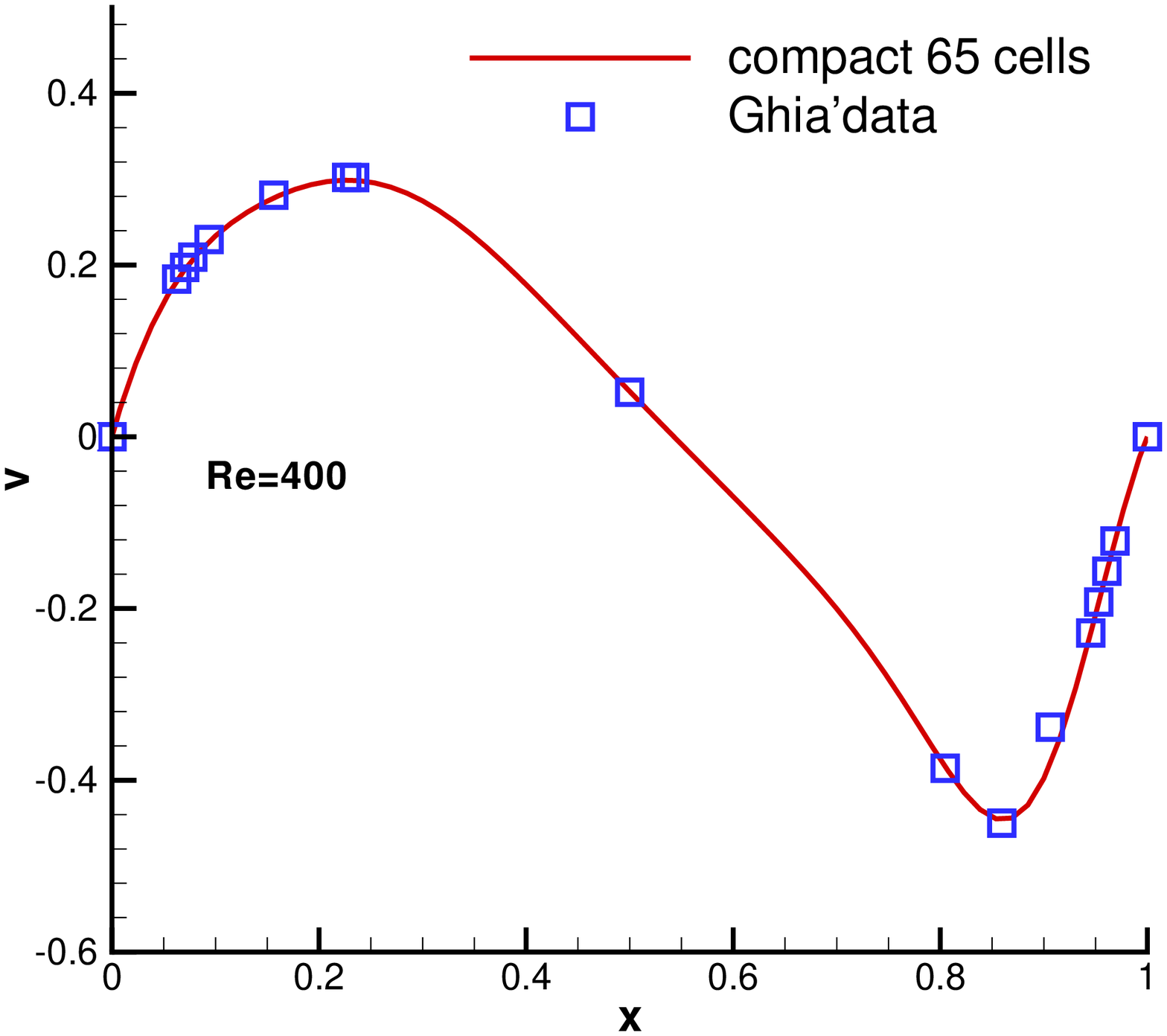}
\caption{\label{cavity-2} Lid-driven cavity flow: $U$-velocities
along vertical centerline line (left) and $V$-velocities along
horizontal centerline with $Re=3200$ and $400$ with $65\times65$
mesh points. The reference data is from \cite{Case-Ghia}.}
\end{figure}

The lid-driven cavity problem is one of the most important
benchmarks for validating incompressible or low speed Navier-Stokes
flow solvers. The fluid is bounded by a unit square and driven by a
uniform translation of the top boundary. In this case, the flow is
simulated with Mach number $Ma=0.15$ and $\gamma=5/3$ in the
computational domain $[0, 1]\times[0, 1]$ and all boundaries are
isothermal and nonslip. Numerical simulations are conducted for
three Reynolds numbers $Re=400, 1000$ and $3200$. The streamlines
with $Re=1000$ using $65\times65$ mesh points are shown in
Fig.\ref{cavity-stream}. The results of $U$-velocities along the
center vertical line, $V$-velocities along the center horizontal
line and the benchmark data \cite{Case-Ghia} are shown in
Fig.\ref{cavity-1} and Fig.\ref{cavity-2}  for different Reynolds
numbers. The simulation results match well with the benchmark data.

\begin{figure}[!h]
\centering
\includegraphics[height=0.4\textwidth]{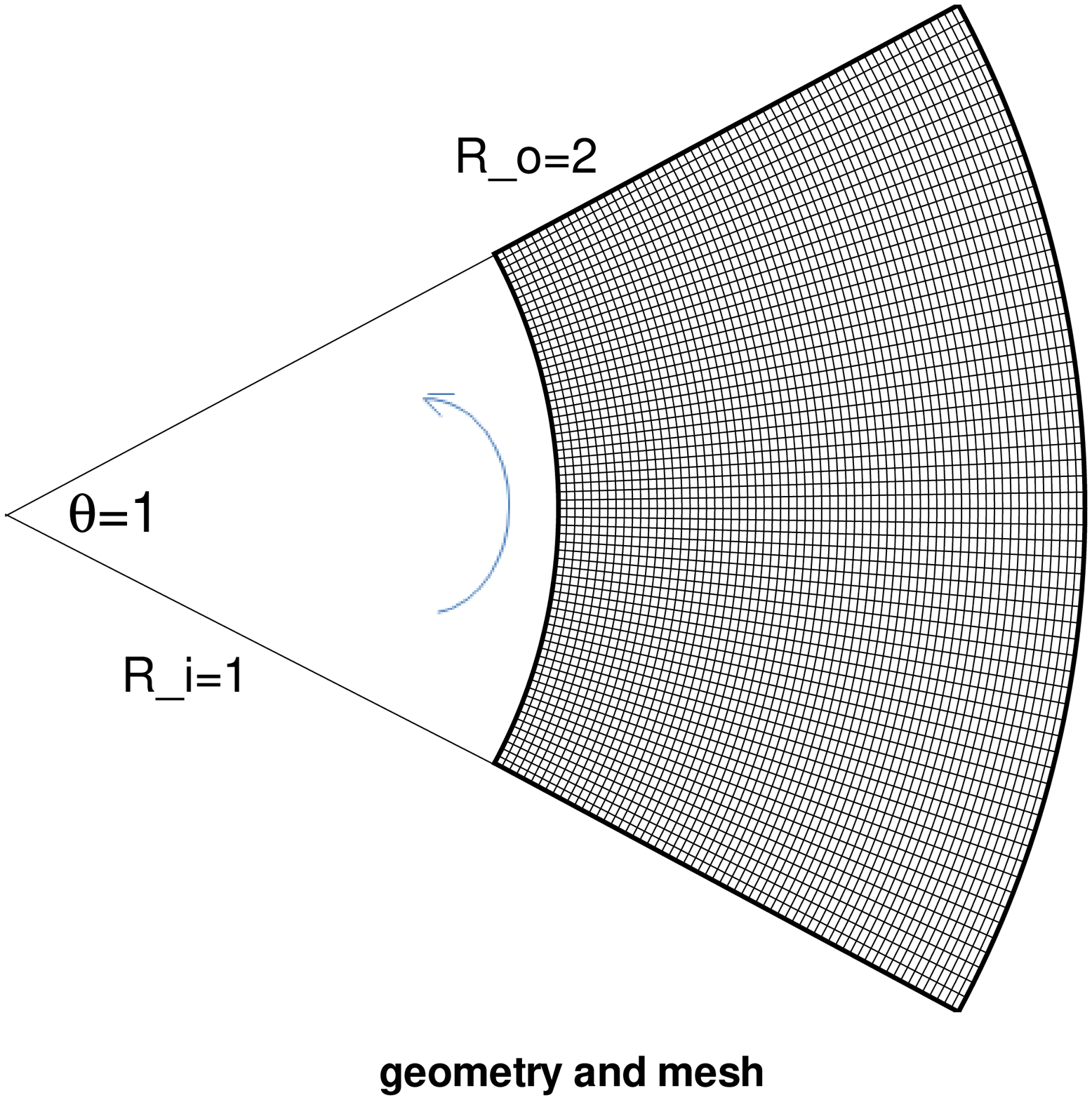}
\includegraphics[height=0.4\textwidth]{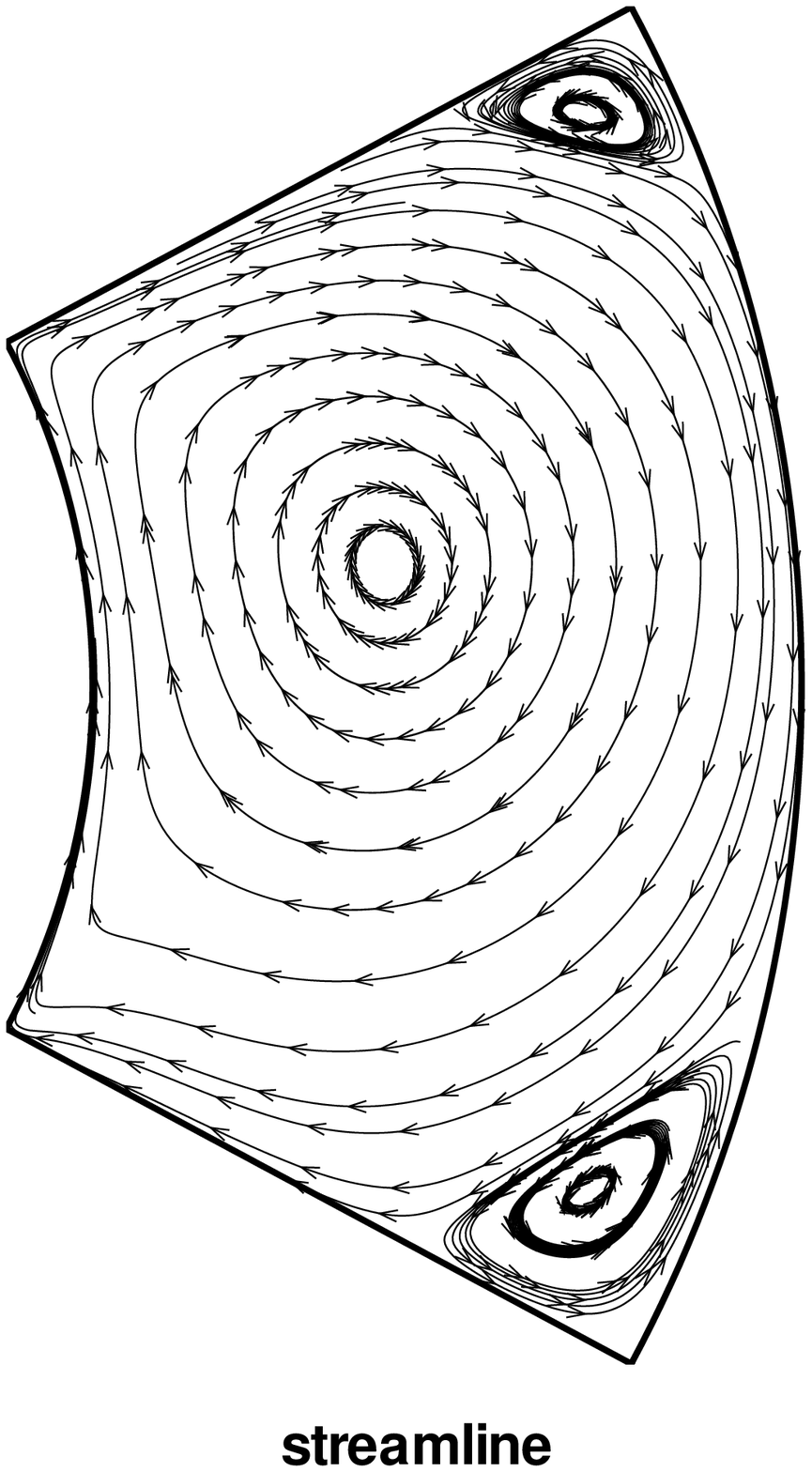}
\caption{\label{cavity-polar1} The schematic diagram (left) and
steady-state streamline profile (right) for the polar cavity flow.}
\end{figure}

\begin{figure}[!h]
\centering
\includegraphics[width=0.42\textwidth]{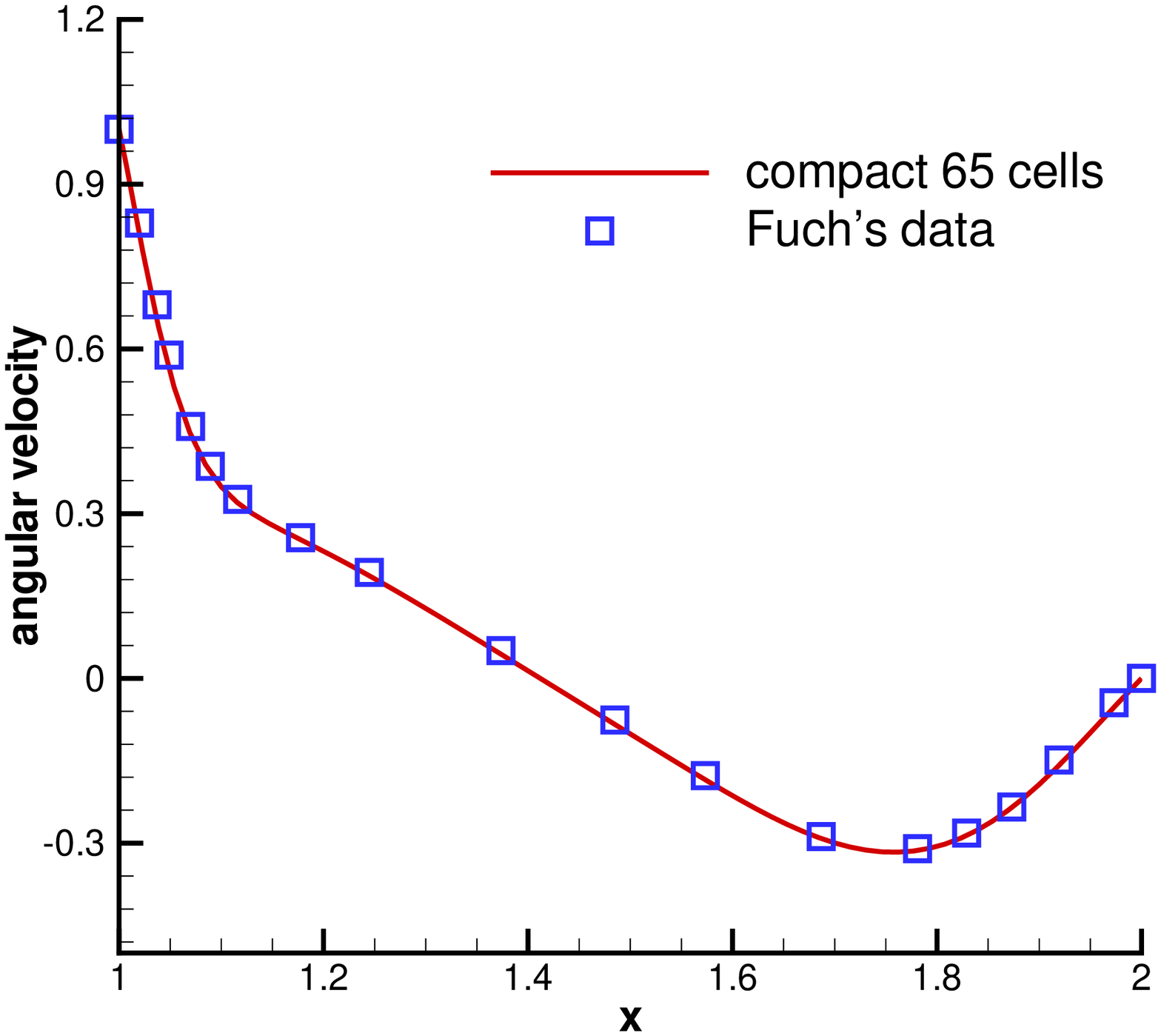}
\includegraphics[width=0.42\textwidth]{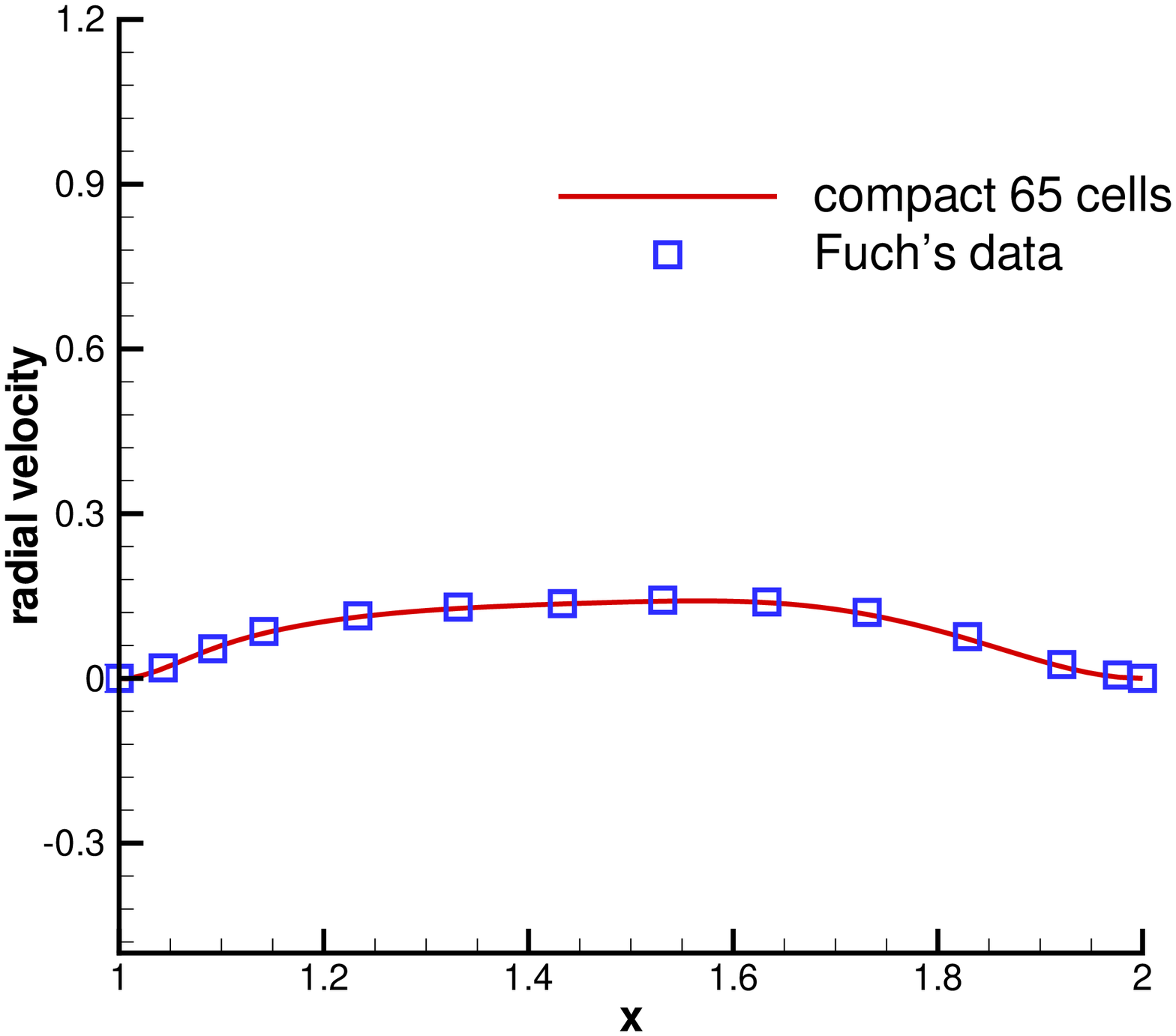}
\caption{\label{cavity-polar2} The profiles of angular velocity and
radical velocity along the horizontal centerline with $\theta=0$ for
the polar cavity flow and the reference data is taken from Fuchs
\cite{Case-Fuchs}.}
\end{figure}

The lid-driven polar cavity flow is tested under the curvilinear
coordinate. The schematic diagram and the computational mesh for
this case are given in Fig.\ref{cavity-polar1}, the computational
domain in the polar coordinate $(r,\theta)$ takes $[1,2]\times[-0.5,
0.5]$, and $65\times65$ uniform mesh points in the polar coordinate
are used. The inner curved wall rotates anticlockwise with
$Ma=0.15$. The flow pattern of this problem is governed by the
Reynolds number defined as $Re =U_{i}R_i/\mu=350$, where $U_i$ is
the inner azimuthal velocity. All the boundaries are also isothermal
and nonslip. The steady-state streamlines are presented in
Fig.\ref{cavity-polar1}. The angular and radial velocity profiles
along the horizontal line with $\theta=0$ are shown in
Fig.\ref{cavity-polar2} together with Fuchs' results
\cite{Case-Fuchs}. Good agreement has been achieved between the
current results and the benchmark data.

\section{Conclusion}
In this paper, a third-order compact gas-kinetic scheme is proposed
for both inviscid and viscous flow computations.
The merit of the current kinetic scheme is that a high-order gas evolution model is constructed and used for
the evaluation of numerical fluxes and the pointwise flow variables at a cell interface.
This can be only achieved with high-order gas evolution model because the solution at the cell interface at the next time
level is a high accurate strong solution of the governing equations. Therefore, with the inclusion of the updated cell interface values,
the scheme can be designed compactly. This can be hardly achieved for the schemes based on the Riemann solution. Therefore, other compact schemes, such as
DG, are constructed based on the weak formulation. Physically, there may have intrinsic inconsistency between the first-order flow dynamics in the Riemann solution
and the high-order flow solver with compact stencil. The weakness in the Riemann dynamics in the traditional ENO and WENO  schemes
are compensated through the large stencils. For example, the high-order derivatives are constructed from the data in the neighboring and neighboring cells in ENO formulation,
instead of updated in the compact DG scheme. Therefore, the compact DG formulation theoretically needs a high-order gas evolution solution, which couples the spatial and
temporal evolution of flow variables compactly.
We believe that the current DG method based on the first-order Riemann solver has intrinsic dynamic weakness in the discontinuous flow regions. In other words,
the weak formulation, which supplies the lost dynamics in the Riemann solver, will be problematic in the regions with discontinuities.
This is probably the reason that the DG can get failed mysteriously for the flow simulation with shocks and complicated flow interactions.
In smooth regions, any governing equation can be manipulated correctly in a physically consistent way, such as all kinds of equivalent weak formulations \cite{CPR3,CPR4}.

 In the current compact gas-kinetic
scheme, both numerical fluxes and pointwise values are used in the construction of the numerical flow evolution.
The core of the scheme is the use of the strong solution of the governing equation from an initial  high-order reconstruction.
Different from the traditional upwind and central schemes, the kinetic formulation is multidimensional, inviscid and viscous terms coupling,
 and has  multi-scale evolution process from the kinetic to the hydrodynamic in the flux construction. This transition in different scale physics is
 equivalent to the transition from the initial upwind scheme to the final central difference one.
 In the current scheme, due to the high-order accuracy in space and time,
the flux transport along a cell interface within a time step can be integrated analytically.
The third-order kinetic scheme doesn't need to use Gaussian points flux integration and the Runge-Kutta time stepping.
This scheme has been validated through the computations for the flows from the smooth incompressible to the hypersonic viscous interaction.
Due to the high-order dynamics in the gas-kinetic formulation, more information can be extracted at the cell interface.
How to use these information is an interesting research topic for the development of high-order schemes.

\section*{Acknowledgement}
The current work was supported by Hong Kong research grant council
(621011, 620813, 16211014) and HKUST (IRS15SC29,SBI14SC11).

\end{document}